\documentclass[a4paper,12pt]{amsart}

\numberwithin{equation}{section}
\theoremstyle{plain}
\newtheorem{theorem}{\sc \bf Theorem}[section]
\newtheorem{lemma}[theorem]{\sc \bf Lemma}
\newtheorem{corollary}[theorem]{\sc \bf Corollary}
\newtheorem{proposition}[theorem]{\sc \bf Proposition}
\newtheorem{claim}{\sc \bf Claim}
\theoremstyle{definition} 
\newtheorem{definition}[theorem]{\sc \bf Definition}
\newtheorem{remark}[theorem]{\sc \bf Remark}
\newtheorem{example}[theorem]{\sc \bf Example}

\renewcommand{\proofname}{\it Proof.} 

\usepackage{amssymb}
\usepackage{amsfonts}
\usepackage{amsmath}
\usepackage{mathrsfs}
\usepackage{latexsym}
\usepackage{accents}
\usepackage{enumerate}
\usepackage{ulem}
\usepackage{fancybox}
\usepackage{wasysym}
\usepackage{color}
\usepackage{graphicx}
\usepackage{cancel}

\usepackage{tikz}

\newif\ifINITPACKNOTE
\INITPACKNOTEtrue

\newif\ifINITPACKKOUGIFUNC
\INITPACKKOUGIFUNCfalse

\makeatletter
\newenvironment{proofex}[1][\proofname]{\par
  \normalfont
  \topsep6\p@\@plus6\p@ \trivlist
  \item[\hskip\labelsep{\bfseries #1}\@addpunct{\bfseries}]\ignorespaces
}{
   \popQED\endtrivlist\@endpefalse
}
\makeatother

\newcommand{\K}{{\mathbb{K}}}

\let\R=\relax
\let\C=\relax
\let\Z=\relax

\newcommand{\R}{{\mathbb{R}}}
\newcommand{\C}{{\mathbb{C}}}
\newcommand{\Z}{{\mathbb{Z}}}

\definecolor{COLgrey}{rgb}{0.65, 0.65, 0.65}
\definecolor{COLblue}{rgb}{0.65, 0.25, 0.25}
\definecolor{brickred}{cmyk}{0,0.89,0.94,0.28}

\definecolor{skyblue}{rgb}{0.82,0.94,1.0} 
\definecolor{blackblue}{rgb}{0.2,0.2,0.6}

\newcommand{\Si}{\Sigma}

\newcommand{\si}{\sigma}

\newcommand{\ph}{\varphi}

\newcommand{\phe}{\varphi(\epsilon,\cdot)}

\ifINITPACKNOTE
\let\lam=\relax
\let\e=\relax
\newcommand{\lam}{\lambda}
\newcommand{\e}{\epsilon}
\else
\fi

\let\p=\relax
\newcommand{\p}{{\prime}}
\newcommand{\pp}{{\prime\prime}}

\newcommand{\LR}{\mathcal{L}}

\newcommand{\OR}{\mathcal{O}}
\newcommand{\PR}{\mathcal{P}}

\newcommand{\RR}{\mathcal{R}}

\newcommand{\ali}[1]{\begin{align*}#1\end{align*}}
\newcommand{\alil}[1]{\begin{align}#1\end{align}}
\newcommand{\ite}[1]{\begin{enumerate}[(1)]#1\end{enumerate}}
\newcommand{\ites}{\begin{enumerate}}
\newcommand{\itee}{\end{enumerate}}
\let\prop=\relax
\let\cor=\relax
\newcommand{\prop}[1]{\begin{proposition}#1\end{proposition}}
\newcommand{\props}{\begin{proposition}}
\newcommand{\prope}{\end{proposition}}
\newcommand{\cor}[1]{\begin{corollary}#1\end{corollary}}
\newcommand{\cors}{\begin{corollary}}
\newcommand{\core}{\end{corollary}}
\newcommand{\thm}[1]{\begin{theorem}#1\end{theorem}}
\newcommand{\thms}{\begin{theorem}}
\newcommand{\thme}{\end{theorem}}
\newcommand{\lem}[1]{\begin{lemma}#1\end{lemma}}
\newcommand{\lems}{\begin{lemma}}
\newcommand{\leme}{\end{lemma}}

\newcommand{\defis}{\begin{definition}}
\newcommand{\defie}{\end{definition}}

\newcommand{\exams}{\begin{example}}
\newcommand{\exame}{\end{example}}
\newcommand{\rem}[1]{\begin{remark}\normalfont #1\end{remark}}

\newcommand{\pros}{\begin{proof}}
\newcommand{\proe}{\end{proof}}
\newcommand{\pross}{\begin{proofex}}
\newcommand{\prossb}{\begin{proofex}\bf}
\newcommand{\smilemk}{\ifmmode\else\leavevmode\unskip\penalty9999\hbox{}\nobreak\hfill\fi\quad\hbox{\smileyex}}
\newcommand{\proes}{\smilemk\end{proofex}}
\newcommand{\prossn}{\begin{proofex}}
\newcommand{\proesn}{\end{proofex}}

\newcommand{\prossq}{\begin{proof}
\begin{enumerate}
\def\theenumi{(\arabic{enumi})}
\def\labelenumi{\theenumi}
\setlength{\leftmargin}{5pt}	
\setlength{\itemsep}{3pt}		
\setlength{\parskip}{0.0pt}		
\setlength{\itemindent}{15pt}	
\setlength{\leftskip}{-35pt}		
\setlength{\labelsep}{3pt}		
}
\newcommand{\case}[1]{\begin{cases}#1\end{cases}}
\newcommand{\cd}{\cdot}
\newcommand{\up}{\upsilon}
\newcommand{\om}{\omega}

\newcommand{\ti}[1]{\tilde{#1}}

\newcommand{\supp}{\mathrm{supp}}
\newcommand{\di}{\displaystyle}

\newcommand{\dist}{{\mathrm{dist}}}
\newcommand{\diam}{{\mathrm{diam}}}
\newcommand{\var}{\mathrm{var}}

\let\MARU=\relax
\newcommand{\MARU}[1]{{\ooalign{\hfil#1\/\hfil\crcr\raise.167ex\hbox{\mathhexbox20D}}}}

\let\to=\relax
\newcommand{\to}{\ \rightarrow\ }

\newcommand{\qqquad}{\quad\qquad}

\newcounter{constants}
\makeatletter
\def\addconst{
\addtocounter{constants}{1}			
\def\@currentlabel{\arabic{constants}}	
\@currentlabel							
}
\makeatother
\newcommand{\adl}[1]{\addconst\label{c:#1}}
\newcommand{\adr}[1]{\ref{c:#1}}

\begin{document}
\title[Perturbed infinite graph-directed IFS]{Degenerate perturbations of infinite graph-directed iterated function systems}
\address{
{\rm Haruyoshi Tanaka}\\
Course of Mathematics Education\\
Education Graduate School of Education, Naruto University of Education\\
748, Nakajima, Takashima, Naruto-cho, Naruto-City, Tokushima, 772-8502, Japan
}
\email{htanaka@naruto-u.ac.jp}
\author{Haruyoshi Tanaka}
\begin{abstract}
We study infinite graph-directed iterated function systems (GIFS) whose underlying graph is not strongly connected and has countably many vertices and edges. In addition to a summability condition for the physical potential, we provide lower and upper estimates of the Hausdorff dimension of the limit set of such GIFS. Bowen type formula is also given under conformal condition and suitable separate conditions. We also introduce perturbed GIFS in which the images of arbitrarily chosen contraction mappings shrink to a single point. In other words, the graph of the perturbed GIFS differs from that of unperturbed GIFS. Assuming suitable continuity condition on contraction mappings, we prove that the Hausdorff dimension of the limit set of the perturbed GIFS converges to that of the unperturbed GIFS, This result generalizes for finite graphs in \cite{T2019,T2016} to the infinite graph setting. As applications, we consider a perturbed nonconformal mapping, as well as convergence and non-convergence in the Hausdorff dimension for perturbed complex continued fractions with degeneration.
\end{abstract}
\maketitle
\setcounter{tocdepth}{3}
\section{Introduction and outline of our results}\label{sec:intro}
We study the Hausdorff dimension of the limit set of graph-directed iterated function systems (GIFS for short), which were first treated as graph directed Markov systems in \cite{MU}. In particular, we consider graphs whose transition matrices may not be irreducible (see Section \ref{sec:GIFS} for details), and whose vertex and edge sets are countably infinite.
\smallskip
\par
We present the following three main results (1)-(3):
\ite
{
\item[(1)] Under no assumptions on the graph and under mild conditions on nonconformal contraction mappings, we provide an upper estimate for the Hausdorff dimension of the limit set (Theorem\ref{th:GIFS_upper_anygraph}). In addition, under the assumption that the potential is summable and the strong separation condition (SSC) holds, we also give a lower estimate (Theorem\ref{th:GIFS_lower}).
\item[(2)] Without any assumptions on the graph, and assuming that the mappings are conformal and satisfies a separation condition (either SSC or OSC), we establish a Bowen-type formula (Theorem \ref{th:GIFS_Bowen}).
\item[(3)] We formulate a perturbation of GIFS with degeneration, meaning that the images of certain contraction mappings shrink to a single point, resulting in the graph of a perturbed system differs from that of the unperturbed system. We provide sufficient conditions for the continuity of the dimension under degeneration (Theorem \ref{th:conv_dim}). 
}
As applications, we present a concrete example involving a non-finitely irreducible transition matrix (Section \ref{sec:ex_non-finitely}), an example of perturbed affine maps with degeneration (Section \ref{sec:PAMD}), and a study of the continuity of the Hausdorff dimension for perturbed complex continued fractions with degeneration (Section \ref{sec:PCCF}). Moreover, we give an example in which the Hausdorff dimension is discontinuous under a similar perturbed setting (Proposition \ref{prop:per_CF_dege_discon}).
\smallskip
\par
Iterated function systems equipped with directed graphs have been studied extensively, particularly with respect to the Hausdorff dimension of their limit sets. Mauldin and Urba\'nski \cite{MU} considered the case where the edge set is infinite and the vertex set is finite, so that the transition matrix of the graph is finitely irreducible. In contrast, \cite{SU} investigated systems with an infinite vertex set and a finitely irreducible transition matrix. The case of finite graphs that are not strongly connected has also been studied by various authors (e.g., \cite{GMR}). In this paper, we treat systems with countably many vertices and edges without any assumptions of the graph.
A formulation of perturbed GIFS with degeneration was first introduced in \cite{MT2008}, where a cookie-cutter type degenerate perturbation was considered. Degenerate perturbations in finite graph GIFS were also studied in \cite{T2016}. We extend these results to the case of infinite graphs. In particular, we prove the continuity of the Hausdorff dimension using our main result (Theorem \ref{th:GIFS_Bowen}) along with a result concerning transfer operators associated with topological Markov shifts with holes (Theorem \ref{th:conv_dim}). This continuity result serves as a foundation for future work on the convergence and convexity of perturbed equilibrium states splitting into a countable number of equilibrium states. 
\smallskip
\par
To prove the Bowen formula, the assumption of finite irreducibility is often used to approximate the pressure of a potential by restricting it to a finite subgraph (see, for example, \cite[Theorem 2.1.5]{MU}). This approximation plays an important role in the lower bounds of the Hausdorff dimension. In order to consider the non-finitely irreducible case, we develop a method based on transfer operators under a summability condition for the physical potential, building on previous work in \cite{T2024_sum}. This allows us to derive a formula for approximating the pressure (Theorem \ref{th:P(ph)=supP(phn)}). Moreover, we will describe that any GIFS with multi graph is reduced to GIFS with simple graph (Section \ref{sec:reduce_ms}). Simple graph has the advantage that, due to the simplicity of its transition matrix, the topological Markov shift generated by its strongly connected components becomes easier to understand (see (\ref{eq:AH=MT})).
Using these results, we prove our main results.
\smallskip
\par
In Section \ref{sec:GIFS}, we introduce the notation for graph-directed iterated function systems. In Section \ref{sec:reduce_ms}, we describe a method for reducing a GIFS with a multigraph to one with a simple graph.
Section \ref{sec:trans} presents the transfer operator and the associated thermodynamic formalism, which are used to prove our main results. We also formulate a perturbed GIFS with degeneration in Section \ref{sec:GIFSdege}.
Section \ref{sec:dimfor} is devoted to the proofs of our main results. In Section \ref{sec:app}, we discuss a concrete example involving a non-finitely irreducible transition matrix (Section \ref{sec:ex_non-finitely}), a perturbed system of affine maps with degeneration (Section \ref{sec:PAMD}), and a perturbed system of complex continued fractions with degeneration (Section \ref{sec:PCCF}). 
\smallskip
\\
\noindent
{\it Acknowledgment.}\ 
This study was supported by JSPS KAKENHI Grant Number 23K03134. 
\section{Preliminaries}\label{sec:pre}
\subsection{Graph-directed iterated function systems}\label{sec:GIFS}
Let $D$ be a positive integer. Consider a system $(G,(J_{v}),(O_{v}),(T_{e}))$ satisfying the following conditions (G.1)-(G.4):
\begin{itemize}
\item[(G.1)] $G=(V,E,i(\cd),t(\cd))$ is a directed multigraph endowed with countable vertex set $V$, countable edge set $E$, and two maps $i(\cd)$ and $t(\cd)$ from $E$ to $V$. For each $e\in E$, $i(e)$ is called the initial vertex of $e$, and $t(e)$ is called the terminal vertex of $e$.
\item[(G.2)] For each $v\in V$, $J_{v}$ is a compact and connected subset of $\R^{D}$ satisfying that it is equal to the closure of its interior $\overline{\mathrm{int}(J_{v})}$, and $c_{\adl{supJ}}:=\sup_{v\in V}\diam(J_{v})<\infty$.
\item[(G.3)] For each $v\in V$, $O_{v}$ is a bounded, open, and connected subset of $\R^{D}$ such that $c_{\adl{cJO}}:=\inf_{v\in V}\dist(J_{v},\R^{D}\setminus O_{v})>0$.
\item[(G.4)] For each $e\in E$, $T_{e}$ is a $C^{1}$-diffeomorphism map from $O_{t(e)}$ to an open subset of $O_{i(e)}$ with $T_{e}(J_{t(e)})\subset J_{i(e)}$. Moreover, there exists $r\in (0,1)$ such that $\sup_{e\in E}\sup_{x\in O_{t(e)}}\|T_{e}^\p(x)\|\leq r$, where $\|T_{e}^\p(x)\|$ means the operator norm of $T_{e}^\p(x)$. 
\end{itemize}
We refer to a system $(G,(J_{v}),(O_{v}),(T_{e}))$ satisfying (G.1)-(G.4) as a {\it graph iterated function system} ({\it GIFS}). If $G$ is simple, that is, for any $v,u\in V$, there is at most one edge $e$ such that $i(e)=v$ and $t(e)=u$, then we call the system a {\it simple GIFS}. In contrast, when $G$ is multigraph, we call the system a {\it multi GIFS}. We will show in Section \ref{sec:reduce_ms} that any multi GIFS can be reduced to a simple GIFS.

First, we present a coding of multi GIFS. The {\it transition matrix} $A=A^{E}$ associated with $E$ is a zero-one matrix indexed by $E$ such that $A(ee^\p)=1$ if $t(e)=i(e^\p)$ and $A(ee^\p)=0$ if $t(e)\neq i(e^\p)$. The associated code space is defined by
\alil
{
X^{E}=\{\om\in {\textstyle\prod_{n=0}^{\infty}E}\,:\,A^{E}(\om_{n}\om_{n+1})=1 \text{ for all }n\geq 0\}.\label{eq:E^inf=}
}
We often write $\om\in X^{E}$ as $\om=\om_{0}\om_{1}\cdots$. Note that this set is called the topological Markov shift (TMS) with countable state space $E$ and with transition matrix $A^{E}$.

Next, we describe a coding of simple GIFS. In this case, the edge set $E$ is regarded as a subset of $V^{2}$, and we may write $T_{e}$ as $T_{vv^\p}$ if $e=vv^\p\in E$. The transition matrix $A=A^{V}$ based on $V$ is a zero-one matrix indexed by $V$ such that $A(vv^\p)=1$ if $vv^\p\in E$. The corresponding code space is defined by
\alil
{
X^{V}=\{\om\in {\textstyle\prod_{n=0}^{\infty}V}\,:\,A^{V}(\om_{n}\om_{n+1})=1 \text{ for all }n\geq 0\}.\label{eq:V^inf=}
}
The coding map $\pi^{E}\,:\,X^{E}\to \R^{D}$ is well defined by
$\{\pi^{E}\om\}=\bigcap_{k=0}^{\infty}T_{\om_{0}}\circ \cdots\circ \om_{k}(J_{t(\om_{k})})$.
The {\it limit set} of the system is given by the image $K^{E}:=\pi^{E}(X^{E})$. We are interested in the Hausdorff dimension of $K^{E}$. By a similar argument above, if $G$ is simple then the coding $\pi^{V}\,:\,X^{V}\to \R^{D}$ is given by $\{\pi^{V}\om\}=\bigcap_{k=0}^{\infty}T_{\om_{0}\om_{1}}\circ \cdots T_{\om_{k}\om_{k+1}}(J_{\om_{k+1}})$ and let $K^{V}:=\pi^{V}(X^{V})$. Then we notice the equation $K_{V}=K_{E}=:K$.
We put
\alil
{
\ph^{E}(\om)=&\log \|T_{\om_{0}}^{\p}(\pi^{E}\si\om)\|\text{ for }\om\in X^{E},\ \ph^{V}(\om)=\log \|T_{\om_{0}\om_{1}}^{\p}(\pi^{V}\si\om)\|\text{ for }\om\in X^{V}\label{eq:ophiV}\\
\underline{\ph}^{E}(\om)=&\log \|T_{\om_{0}}^{\p}(\pi^{E}\si\om)\|_{i}\text{ for }\om\in X^{E},\ \underline{\ph}^{V}(\om)=\log \|T_{\om_{0}\om_{1}}^{\p}(\pi^{V}\si\om)\|_{i}\text{ for }\om\in X^{V},\label{eq:uophiV}
}
where $\si$ is the shift transformation on $X^{E}$ or $X^{V}$ and we put $\|T^\p_{e}(x)\|_{i}:=\|T^\p_{e}(x)^{-1}\|^{-1}=\inf_{x\in \R^{D}\,:\,|x|\geq 1}|T^\p_{e}(x)|$. We impose the following conditions as necessary.
\ite
{
\item[$(G.5)_{C}$] (Conformality) Each $T_{e}$ is conformal.
\item[$(G.6)_{S}$] (Strong Separation Condition (SSC)) For $e,e^\p\in E$ with $e\neq e^\p$, $T_{e}(J_{t(e)}\cap K)\cap T_{e^\p}(J_{t(e^\p)}\cap K)=\emptyset$.
\item[$(G.6)_{O}$] (Open Set Condition (OSC)) For $e,e^\p\in E$ with $e\neq e^\p$, $T_{e}(\mathrm{int}J_{t(e)})\cap T_{e^\p}(\mathrm{int}J_{t(e^\p)})=\emptyset$.
\item[$(G.7)_{B}$] (Bounded Distortion) There exist constants $c_{\adl{bd1}}>0$ and $0<\beta\leq 1$ such that for any $e\in E$ and $x,y\in O_{t(e)}$, $|\|T_{e}^\p(x)\|-\|T_{e}^\p(y)\||\leq c_{\adr{bd1}}\|T_{e}^\p(x)\| |x-y|^{\beta}$.
\item[$(G.7)_{BI}$] (Bounded Distortion for Inverse Maps) There exist constants $c_{\adl{bd2}}>0$ and $0<\beta\leq 1$ such that for any $e\in E$ and $x,y\in O_{t(e)}$, $|\|T_{e}^\p(x)\|_{i}-\|T_{e}^\p(y)\|_{i}|\leq c_{\adr{bd2}}\|T_{e}^\p(x)\|_{i} |x-y|^{\beta}$. 
\item[$(G.8)_{SE}$] (Summability using $E$)
Define
\ali
{
\underline{I}^{E}:=&\{s\geq 0\,:\,\sum_{e\in E}\sup_{x\in K\cap J_{t(e)}}\|T_{e}^\p(x)\|_{i}^{s}<\infty\}.
}
Then there exists $s\in \underline{I}^{E}$ such that $P(s\underline{\ph}^{E})=0$, where $P(s\underline{\ph}^{E})$ denotes the topological pressure introduced in (\ref{eq:P(f)=}).
\item[$(G.8)_{SV}$] (Summability using $V$)
If $G$ is simple, define
\ali
{
\underline{I}^{V}:=&\{s\geq 0\,:\,\sum_{v\in V}\sup_{u\in V\,:\,vu\in E}\sup_{x\in K\cap J_{u}}\|T_{vu}^\p(x)\|_{i}^{s}<\infty\}.
}
Then there exists $s\in \underline{I}^{V}$ such that $P(s\underline{\ph}^{V})=0$.
\item[$(G.9)_{J}$] (Contractibility of $J_{v}$) There exists a constant $c_{\adl{CJ}}>0$ such that for each $v\in V$, $\diam J_{v}\leq c_{\adr{CJ}}\sup_{e\in E\,:\,i(e)=v}\sup_{z\in J_{t(e)}}\|T_{e}^\p(z)\|$.
}
The GIFS is said to be {\it strongly regular} if $0<P(s\ph^{E})<+\infty$ or $0<P(s\ph^{V})<+\infty$ for some $s\geq 0$. It is called {\it regular} if $P(s\ph^{E})=0$ or $P(s\ph^{V})=0$ for some $s\geq 0$. Otherwise, it is called {\it irregular}. A subgraph $G_{0}=(V_{0},E_{0},i(\cd),t(\cd))$ of $G$ is {\it strongly connected} if for any two vertexes $v,u\in V_{0}$ there exists a path $w\in \bigcup_{n=1}^{\infty}(E_{0})^{n}$ on $G$ such that $v=i(w)$ and $u=e(w)$. A strongly connected subgraph $G_{0}$ is {\it strongly connected component} if for any strongly connected subgraph $G_{1}$ satisfying that $G_{0}$ is a subgraph of $G_{1}$, we have the equality $G_{1}=G_{0}$. We let
\ali
{
\mathrm{SC}(G):=\{G_{0}\,:\,\text{strongly connected component of }G\}.
}
\rem
{\label{eq:GIFS_rem}
\item[(1)] Instead of the condition $\sup_{e\in E}\sup_{x\in O_{t(e)}}\|T_{e}^\p(x)\|<1$ in (G.4), we may assume the following: there exist constants $c_{\adl{cr2}}>0$ and $r\in (0,1)$ such that for any $n\geq 1$ and any word $w=e_{1}\cdots e_{n}\in E^{n}$,  we have $\sup_{x\in O_{t(w)}}\prod_{i=1}^{n}\|T_{e_{i}}^\p(T_{e_{i+1}}\circ\cdots \circ T_{e_{n}}(x))\|\leq c_{\adr{cr2}}r^{n}$. In fact, all the results in Section \ref{sec:GIFSdege} and Section \ref{sec:dimfor} remain valid under this condition. This condition is applicable to the GIFS associated with complex continued fractions (see also Section \ref{sec:PCCF}).
\item[(2)] If condition $(G.5)_{C}$ is satisfied, then $\|T_{e}^\p(x)\|=\|T_{e}^\p(x)\|_{i}$. Hence, conditions $(G.7)_{B}$ and $(G.7)_{BI}$ become equivalent. 
\item[(3)] If the graph $G$ is simple, then condition $(G.8)_{SE}$ is stronger than condition $(G.8)_{SV}$ (see also Section \ref{sec:ex_non-finitely}).
\item[(4)] Assume that conditions (G.1)-(G.4) and $(G.7)_{B}$ are satisfied. It is known from \cite{MU} that if $A^{E}$ is finitely irreducible, then the equation $\inf \underline{I}^{E}=\inf\{s\geq 0\,:\,P(s\underline{\ph}^{E})<\infty\}$ holds.
Then in this case, the assumption $s\in \underline{I}^{E}$ in condition $(G.8)_{SE}$ is not required. 
\item[(5)] We comment on condition $(G.9)_{J}$. Assume that $G$ is simple. Since the inclusion $\bigcup_{u\in V}T_{vu}J_{u}\subset J_{v}$ and $J_{v}$ is compact, it is natural to expect that if a separation condition holds, then $\diam(T_{vu}J_u)$ must vanish as $u$ becomes large. On the other hand, in our setting, we allow $|T_{vu}^\p|$ to remain uniformly positive with respect to $u\in V$. Therefore, it is necessary that $\diam J_u$ vanishes as $u$ to large, and condition $(G.9)_{J}$ imposes a constraint on the speed of this convergence.
}
\subsection{Reduction from multi GIFS to simple GIFS}\label{sec:reduce_ms}
In this section, we demonstrate that a multi GIFS can be reduced to a simple GIFS. We begin with auxiliary results which are useful for proving our main theorem.
\begin{proposition}[\cite{Patzschke}]\label{prop:T:mean}
Let $(G,(J_{v}),(O_{v}))$ be a system with conditions (G.1)-(G.4). Put $U_{v}=B(J_{v}),c_{\adr{cJO}})$ for $v\in V$. Then there exists a constant $c_{\adl{MT}}\geq 1$ such that for any $v\in V$ and any $C^{1}$ map $T\,:\,O_{v}\to Y$, where $(Y,\|\cdot\|_{Y})$ is a normed space, the following inequality holds for each $x,y\in J_{v}$:
\begin{equation}\label{eq:proof_1}
\|T(x)-T(y)\|_{Y}\leq c_{\adr{MT}}\sup_{z\in U_{v}}\|T^{\prime}(z)\| |x-y|.
\end{equation}
The constant $c_{\adr{MT}}$ depends only on $D, c_{\adr{supJ}}$ and $c_{\adr{cJO}}$. 
\end{proposition}
\prop
{[\cite{Patzschke}]\label{prop:T:mean_inf}
Let $(G,(J_{v}),(O_{v}))$ be a system satisfying conditions (G.1)-(G.3). Then there exists a constant $0<c_{\adl{iMT}}\leq 1$ such that for any $v,v^\p\in V$ and any $C^{1}$-diffeomorphism $T\,:\,O_{v}\to O_{v^\p}$, the inequality $|T(x)-T(y)|\geq c_{\adr{iMT}}\inf_{z\in O_{v}}\|T^\p(z)^{-1}\|^{-1}|x-y|$ holds for any $x,y\in J_{v}$, where the constant $c_{\adr{iMT}}$ depends only on $D, c_{\adr{supJ}}$ and $c_{\adr{cJO}}$.
}
Let $\mathcal{G}=(G,(J_{v}),(O_{v}), (T_{e}))$ be a multi GIFS satisfying (G.1)-(G.4) and either $(G.6)_{O}$ or $(G.6)_{S}$. We construct a simple GIFS $\tilde{\mathcal{G}}=(\tilde{G}=(\tilde{V},\tilde{E}),(\tilde{J}_{v}),(\tilde{O}_{v}),(T_{vv^\p})_{vv^\p\in E})$ as follows:
\ali
{
&\tilde{V}:=E,\ \tilde{E}:=\{ee^\p\in E^{2}\,:\,t(e)=i(e^\p)\}\\
&\tilde{J}_{e}:=T_{e}(J_{t(e)}),\ \tilde{O}_{e}:=O_{i(e)}\quad (e\in \tilde{V})\quad\text{ and }\tilde{T}_{ee^\p}:=T_{e}\quad (ee^\p\in \tilde{E}).
}
\prop
{\label{prop:reduce_GIFS}
Under the above notation, wthe following statements hold:
\ite
{
\item The new system $\tilde{\mathcal{G}}$ is a simple GIFS.
\item $X^{E}=X^{\tilde{V}}$, $\pi^{E}=\pi^{\tilde{V}}$ and $\ph^{E}=\ph^{\tilde{V}}$.
\item If the original system $\mathcal{G}$ satisfies $(G.5)_{C}$, $(G.6)_{S}$, $(G.6)_{O}$, $(G.7)_{B}$ or $(G.7)_{BI}$, then so does the new system $\tilde{\mathcal{G}}$.
\item The interval $\underline{I}^{E}$ of $\mathcal{G}$ is equal to the interval $\underline{I}^{\tilde{V}}$ of the new system $\tilde{\mathcal{G}}$. In particular, if $\mathcal{G}$ satisfies $(G.8)_{SE}$, then $\tilde{\mathcal{G}}$ satisfies $(G.8)_{SV}$.
\item Assume that for any $e\in E$, there exists $e^\p\in E$ such that $t(e)=i(e^\p)$, and that $\mathcal{G}$ satisfies $(G.7)_{B}$. Then $\tilde{\mathcal{G}}$ satisfies $(G.9)_{J}$.
}
}
\pros
(1) We show that the new system satisfies conditions (G.1)-(G.4).
Since $\tilde{E}\subset \tilde{V}^{2}$, the graph $\tilde{G}$ is simple. Thus, (G.1) is satisfied.
Each $T_{e}$ is diffeomorphism, so each set $\tilde{J}_{e}=T_{e}(J_{t(e)})$ is a nonempty, compact and connected set with $\overline{\mathrm{int} \tilde{J}_{e}}=\tilde{J}_{e}$. Moreover, by Proposition \ref{prop:T:mean}, we have $\sup_{e\in \tilde{V}}\diam \tilde{J}_{e}\leq c_{\adr{MT}}r\sup_{v\in V}\diam J_{v}<+\infty$. Therefore (G.2) is fulfilled.
Since $\tilde{J}_{e}=T_{e}(J_{t(e)})\subset J_{i(e)}$, we obtain $\dist (\tilde{J}_{e},\R^{D}\setminus \tilde{O}_{e})\geq \dist (J_{i(e)}\setminus O_{i(e)})\geq \inf_{v\in V}\dist (J_{v}\setminus O_{v})>0$. Thus (G.3) is valid.
To check (G.4), observe that $\tilde{T}_{ee^\p}\tilde{J}_{e^\p}=T_{e}T_{e^\p}J_{t(e^\p)}\subset T_{e}J_{i(e^\p)}=T_{e}J_{t(e)}=\tilde{J}_{e}$. Condition (G.1)-(G.4) are satisfied.\\
(2)(3) These follow directly from the definition of $\tilde{\mathcal{G}}$.\\
(4) Assume that the original system satisfies $(G.8)_{SE}$. Then
\ali
{
\sum_{e\in \tilde{V}}\sup_{e^\p\in \tilde{V}\,:\,\atop{ee^\p\in \tilde{E}}}\sup_{x\in \tilde{J}\cap \tilde{J}_{e^\p}}\|\tilde{T}_{ee^\p}(x)^\p\|_{i}^{s}=\sum_{e\in E}\sup_{e^\p\in E\,:\,\atop{t(e)=i(e^\p)}}\sup_{x\in T_{e^\p}(J_{t(e^\p)})\cap J}\|T_{e}(x)^\p\|_{i}^{s}. 
}
For any $x\in T_{e^\p}(J_{t(e^\p)})$ with $e^\p\in E$ and $t(e)=i(e^\p)$, and for any $y\in J_{t(e)}$, we have $c^{-1}\|T_{e}^\p(y)\|\leq \|T_{e}^\p(x)\|\leq c\|T_{e}^\p(y)\|$ for some $c\geq 1$.
Hence we have $\underline{I}^{E}=\underline{I}^{V}$, and so $(G.8)_{SV}$ holds for the new system.
\smallskip
\\
(5) For each $e\in E$, take $e^\p\in E$ so that $ee^\p\in \tilde{E}$. Then
\ali
{
\diam(\tilde{J}_{e})=\diam(T_{e}(J_{t(e)}))\leq c_{\adr{MT}}\sup_{z\in O_{t(e)}}\|T_{e}^\p(z)\|\diam J_{t(e)}
&\leq c_{\adr{MT}}c_{\adr{supJ}}\sup_{z\in \tilde{O}_{e^\p}}\|T_{ee^\p}^\p(z)\|.
}
Thus, the condition $(G.9)_{J}$ is satisfied for $\tilde{\mathcal{G}}$.
\proe
\subsection{Transfer operators and thermodynamic formalism}\label{sec:trans}
We recall some properties of transfer operators associated with holes, as well as aspects of thermodynamic formalism discussed in \cite{T2024_sum}.
Let $X$ be a topological Markov shift with state space $S$ and transition matrix $A=(A(ij))$ indexed by $S$. The shift transformation $\si\,:\,X\to X$ is defined by $(\si\om)_{n}=\om_{n+1}$ for any $n\geq 0$. Let $S_{0}\subset S$ be a nonempty subset, and let $M=(M(ij))$ be a zero-one matrix indexed by $S_{0}$, such that $M(ij)\leq A(ij)$ for any $i,j\in S_{0}$. Define
\alil
{
\textstyle{X_{M}:=\{\om\in \prod_{n=0}^{\infty}S_{0}\,:\,M(\om_{n}\om_{n+1})=1\text{ for all }n\geq 0\}}.\label{eq:XM=...}
}
In other words, $X_{M}$ is a subsystem of $X$ with state space $S_{0}$ and transition matrix $M$. For convenience, we define $M(ij):=0$ whenever $i\in S\setminus S_{0}$ or $j\in S\setminus S_{0}$.
A word $w=w_{1}w_{2}\dots w_{n}\in S^{n}$ is called {\it $M$-admissible} if $M(w_{i}w_{i+1})=1$ for all $1\leq i<n$.
For $a,b\in S$, we write $a\to b$ if there exist an integer $n\geq 2$ and an $M$-admissible word $w\in S^{n}$ such that $w_{1}=a$ and $w_{n}=b$.
The matrix $M$ is said to be {\it irreducible} if $a\to b$ for any $a,b\in S$. Next, we introduce tools to handle non-irreducible transition matrices. For $a,b\in S$, we write $a\leftrightarrow b$ if either $a\to b$ and $b\to a$, or $a=b$. Then $\leftrightarrow$ is an equivalent relation on $S$. Consider the quotient space $S/\!\!\leftrightarrow$. For $S_{1},S_{2}\in S/\!\!\leftrightarrow$, we write $S_{1}\preceq_{M} S_{2}$ if either $S_{1}=S_{2}$, or there exist $a\in S_{1}$ and $b\in S_{2}$ such that $a\to b$. This defines a semi-order on $S/\!\!\leftrightarrow$. For $T\in S/\!\!\leftrightarrow$, denote by $M(T)$ the irreducible submatrix of $M$ indexed by $T$. Consequently, we obtain countably many transitive components $X_{M(T)}$ of $X$, one for each $T\in S/\!\!\leftrightarrow$.

Assume that $G$ is simple graph. We now describe a relationship between $S/\!\!\leftrightarrow$ and $\mathrm{SC}(G)$. For each $H\in \mathrm{SC}(G)$, let $X_{H}$ be the topological Markov shift (TMS) defined by the transition matrix of $H$ using its vertex set as the state space.
Then we have the following identity:
\alil
{
\{X_{H}\,:\,H\in \mathrm{SC}(G)\}=\{X_{M(T)}\,:\,T\in S/\!\!\leftrightarrow \text{ and }M(T)\text{ is not }1\times 1 \text{ zero matrix}\}.\label{eq:AH=MT}
}
The matrix $A$ is said to be {\it finitely irreducible} if there exists a finite subset $F$ of $\bigcup_{n=1}^{\infty}S^{n}$ such that for any $a,b\in S$, there exists $w\in F$ for which $awb$ is $A$-admissible.
\smallskip
\par
 
A function $\ph\,:\,X\to \R$ is called {\it summable} if
\alil
{
\sum_{s\in S\,:\,[s]\neq \emptyset}\exp(\sup_{\om\in [s]}\ph(\om))<\infty, \label{eq:sum}
}
where for a word $w\in S^{n}$, the cylinder set $[w]$ is defined by $[w]=\{\om\in X\,:\,\om_{0}\cdots\om_{n-1}=w\}$.
For $\theta\in (0,1)$, define a metric $d_{\theta}$ on $X$ by $d_{\theta}(\om,\up)=\theta^{\min\{n\geq 0\,:\,\om_{n}\neq \up_{n}\}}$ if $\om\neq \up$ and $d_{\theta}(\om,\up)=0$ if $\om=\up$.
Let $\K=\R$ or $\C$. For a function $f\,:\,X\to \K$ and an integer $k\geq 1$, define
\ali
{
[f]_{k}=\sup\{\var_{n}f/\theta^{n}\,:\,n\geq k\},
}
where $\var_{n}f:=\sup\{|f(\om)-f(\up)|\,:\,\om_{i}=\up_{i}\ \ (i=0,1,\dots, n-1)\}$. Note that $[f]_{k}\geq [f]_{k+1}$. If $[f]_{1}<\infty$, then $f$ is called {\it locally Lipschitz continuous}, and if $[f]_{2}<\infty$, then $f$ is called {\it weak Lipschitz continuous}.
Let $C(X)$ denote the space of all $\C$-valued continuous functions on $X$, and let $F^{k}(X)$ be the space of all $f\in C(X)$ such that $[f]_{k}<\infty$. Define $C_{b}(X)$ to be the Banach space of all bounded functions $f\in C(X)$, equipped with the supremum norm $\|f\|_{\infty}=\sup_{X}|f|$. Similarly, let $F^{k}_{b}(X)$ be the Banach space of all $f\in F^{k}(X)$ such that $f\in C_{b}(X)$, endowed with the norm $\|f\|_{k}=\|f\|_{\infty}+[f]_{k}$.

Next we recall some properties of thermodynamic formalism and transfer operators. 
For a real-valued function $\ph$ on $X$, the {\it topological pressure} $P(\ph)$ is defined by
\alil
{
P(\ph):=\lim_{n\to \infty}\frac{1}{n}\log \sum_{w\in S^{n}\,:\,[w]\neq \emptyset}\exp(\sup_{\om\in [w]}S_{n}\ph(\om)),\label{eq:P(f)=}
}
where we define $S_{n}\ph(\om):=\sum_{i=0}^{n-1}\ph(\si^{i}\om)$. As shown in \cite[Lemma 2.1.2]{MU}, if $A$ is finitely irreducible and $\ph\,:\,X\to \R$ is locally Lipschitz continuous, then $P(\ph)<+\infty$ if and only if $\ph$ is summable.

A $\sigma$-invariant Borel probability measure $\mu$ on $X$ is called a {\it Gibbs measure} for the potential $\varphi$ if there exist constants $c_{\adl{Gibbs1}}>0$ and $P\in \R$ such that for any $\om\in X$ and $n\geq 1$,
\begin{align}
c_{\adr{Gibbs1}}^{-1}\leq \frac{\mu([\om_{0}\om_{1}\cdots \om_{n-1}])}{\exp(-nP+S_{n}\varphi(\om))}\leq c_{\adr{Gibbs1}}.\label{eq:Gibbs}
\end{align}
It is well known that such a measure exists uniquely when $A$ is finitely irreducible and $\varphi$ is locally Lipschitz continuous. In this case, the constant $P$ coincides with $P(\ph)$.

We now introduce a transfer operator associated with an open system of the shift. The operator $\LR_{M}=\LR_{M,\ph}$ associated with $M$ and $\ph$ is formally defined by
\alil
{
\LR_{M} f(\om)=\LR_{M,\ph} f(\om):=\sum_{a\in S\,:\,M(a\om_{0})=1}e^{\ph(a\cd \om)}f(a\cd\om)\label{eq:transfer}
}
for $f\,:\,X\to \C$ and $\om\in X$, where $a \cd \om$ denotes the concatenation $a\om_{0}\om_{1}\cdots$. 
Note that the operator $\LR_{M}$ is bounded acting on $F^{k}_{b}(X)$ and on $C_{b}(X)$. 
Such operators are used in studying system with holes \cite{Demers_,T2024_sum, T2020, T2019, T2009}. In fact, by setting $\Si=\bigcup_{ij\,:\,M(ij)=1}[ij]$, we regard the map $\si|_{\Si}\,:\,\Si\to X$ as an {\it open system}, while the map $\si\,:\,X\to X$ is viewed as a {\it closed system}. 

Assume the following three conditions:
\ite
{
\item[(A.1)] The matrix $A$ is irreducible. 
\item[(A.2)] The subsystem $X_{M}$ of $X_{A}$ with the transition matrix $M$ has a periodic point for $\si$.
\item[(A.3)] A function $\ph\,:\,X\to \R$ is summable and satisfies $[\ph]_{k+1}<+\infty$ for some $k\geq 1$.
}
\props
[{\cite[Proposition 4.2]{T2024_sum}}]\label{prop:positive_eval}
Assume that the conditions (A.1) and (A.3) are satisfied. Then condition (A.2) is satisfied if and only if the spectral radius of $\LR_{M,\ph}\,:\,C_{b}(X)\to C_{b}(X)$ is positive. 
\prope
\thms[\cite{T2024_sum}]
\label{th:ex_efunc_geneRuelleop_M}
Assume that conditions (A.1)-(A.3) are satisfied, and that $M$ is irreducible. Let $S_{0}$ be the index set of $M$. Then the operator $\LR_{M}\,:\,F^{k}_{b}(X)\to F^{k}_{b}(X)$ admits the spectral decomposition
\ali
{
\LR_{M}=\lam\PR+\RR
}
satisfying the following:
\ite
{
\item[(i)] $\lam$ is a positive eigenvalue of $\LR_{M}$, equal to its spectral radius, and given by $\lam=\exp(P(\ph|_{X_{M}}))$.
\item[(ii)] $\PR$ is the projection onto the one-dimensional eigenspace corresponding to $\lam$, and has the form $\PR=h\otimes \nu$, where $h\in F_{b}^{k}(X)$ is a nonnegative eigenfunction supported on $\bigcup_{a\in S_{0}}[a]$, and $\nu$ is the corresponding positive eigenvector of the dual operator $\LR_{M}^{*}$, normalized by $\nu(1)=\nu(h)=1$.
\item[(iii)] The remainder $\RR$ satisfies $\RR\PR=\PR\RR=\OR$, and its spectrum is in the outside of $\{\lam\}$.
}
\thme
\pros
See Theorem 3.1, Proposition 4.6, and Proposition 4.7 in \cite{T2024_sum}.
\proe
For simplicity, we refer to the triple $(\lam,h,\nu)$ as the {\it spectral triplet} of $\LR_{M}$.
\thm
{\label{th:pressure_max}
Let $X$ be a TMS with transition matrix $A$ and state space $S$. Assume that condition (A.3) is satisfied. If $X$ has a periodic point under $\si$, then $P(\ph)=\max_{T\in S/\leftrightarrow}P(\ph|_{X_{M(T)}})$. If $X$ has no periodic points, then $P(\ph)=-\infty$.
}
\pros
Let $\hat{A}$ be an $S\times S$ matrix whose entries are all $1$, and let $\hat{X}$ be the TMS with the transition matrix $\hat{A}$.
By using \cite[Appendix A]{T2024_sum}, there exists a real-valued function $\hat{\ph}\in F_{b}^{k+1}(\hat{X})$ such that $\hat{\ph}=\ph$ on $X$ and $\hat{\ph}$ is summable. First we assume that $X$ has a periodic point under $\si$. Then conditions (A.1)-(A.3) are satisfied replacing $A$ and $M$ by $\hat{A}$ and $A$, respectively. Thus it follows from \cite[Proposition 4.9]{T2024_sum} that $P(\ph)=P(\hat{\ph}|X_{A})=\max_{T\in S/\leftrightarrow}P(\ph|X_{M(T)})$. On the other hand, if there are no periodic points of $X$, then Proposition \ref{prop:positive_eval} implies that $P(\ph)=P(\hat{\ph}|_{X_{A}})\leq \log \lim_{n\to \infty}\|\LR_{A}^{n}\|_{\infty}^{1/n}=-\infty$. Hence the assertion is valid.
\proe
We also recall the following condition for perturbed potentials depending on a small parameter $\e\in (0,1)$ (see Section 4 in \cite{T2024_sum}):
\ite
{
\item[(B.1)] Real-valued functions $\phe\in F^{k+1}(X)$ ($\e\in (0,1)$) satisfy $\sup_{\e>0}[\ph(\e,\cd)]_{k+1}<\infty$.
\item[(B.2)] $c_{\adl{us}}:=\sum_{s\in S}\exp(\sup_{\e>0}\sup_{\om\in [s]}\ph(\e,\om)<\infty$.
\item[(B.3)] Take $\psi\,:\,X\to \R$ so that $\psi(\om)=\exp(\ph(\om))$ if $M(\om_{0}\om_{1})=1$ and $\psi(\om)=0$ otherwise. Then $\sup_{\om\in [a]}|\exp(\ph(\e,\om))-\psi(\om)|\to 0$ as $\e\to 0$ for each $a\in S$.
}
Note that condition (B.1) above is weaker than the corresponding condition in \cite{T2024_sum}, where $\phe$ is required to belong to $F^{1}(X)$. Recall also that there exists a perturbed potential $\phe$ satisfying conditions (B.1)-(B.3) (see \cite[Lemma 5.1]{T2024_sum}).
\thm
{\label{th:perTMS}
Assume that conditions (A.1)-(A.3) and (B.1)-(B.3) are satisfied. Then 
\ite
{
\item $\|\LR_{A,\phe}-\LR_{M}\|_{\infty}\to 0$.
\item Assume further that $M$ is irreducible. Let $(\lam(\e),h(\e,\cd),\nu(\e,\cd))$ be the spectral triplet of $\LR_{A,\phe}$, and $(\lam,h,\nu)$ that of $\LR_{M}$. Then we have convergence $\lam(\e)\to \lam$, $h(\e,\om)\to h(\om)$ for each $\om\in X$, and $\nu(\e,f)\to \nu(f)$ for each $f\in C_{b}(X)$, as $\e\to 0$.
\item $P(\phe)\to P(\ph)$ as $\e\to 0$.
}
}
\pros
(1) See Proposition 4.3 in \cite{T2024_sum}. Although the original result was proved under the assumption $\phe\in F^{1}(X)$, the proof remains valid for $\phe\in F^{k+1}(X)$ without modification. Following (2) and (3) are also similar.\\
(2) See Proposition 4.4 and Proposition 4.6 in \cite{T2024_sum}, together with Theorem \ref{th:ex_efunc_geneRuelleop_M}.\\
(3) See Proposition 6.1 in \cite{T2024_sum}.
\proe

We now present the following useful result for the proof of Theorem \ref{th:P(ph)=supP(phn)}. 
\prop
{\label{prop:inc_irremat}
Let $X$ be a TMS with countable state space $S=\{1,2,\dots\}$ and transition matrix $A$. Assume that $A$ is irreducible. Then for any $n\geq 1$, there exists a finite subset $S_{n}$ with $\{1,2,\dots, n\}\subset S_{n}$ and an irreducible zero-one matrix $A_{n}$ indexed by $S_{n}\times S_{n}$ such that $A_{n}=A$ on $\{1,2,\dots, n\}\times S_{n}$, $S_{n}\subset S_{n+1}$, and $A_{n}(ij)\leq A_{n+1}(ij)\leq A(ij)$ for any $n\geq 1$ and $i,j\in S_{n}$.
}
\pros
Let $T_{n}:=\{1,2,\dots, n\}$ and $B_{n}:=A|_{T_{n}\times T_{n}}$. We will construct a family of words $\{w(n,ij)\}$ inductively as follows.
Since $A$ is irreducible, for any $i,j\in T_{1}$, there exist an integer $m(ij)\geq 0$ and a word $w(ij)\in S^{m(ij)}$ such that the concatenation $i\cd w(ij)\cd j$ is $A$-admissible. Now fix $n\geq 2$ and $i,j\in T_{n}$. If $i\in T_{n}\setminus T_{n-1}$ or $j\in T_{n}\setminus T_{n-1}$, then we take an integer $m(ij)\geq 0$ and a word $w(ij)\in S^{m(ij)}$ such that $i\cd w(ij)\cd j$ is $A$-admissible. 
Put
$S_{n}:=T_{n}\cup\{w(ij)_{1}\,:\,i,j\in T_{n}\}$ and define $A_{n}\,:\,S_{n}\times S_{n}\to \{0,1\}$ by
\ali
{
A_{n}(ij)=
\case
{
A(ij), &\text{if }i\in T_{n}\\
1,&\text{if }i\in S_{n}\setminus T_{n} \text{ and }ij=w(i^\p j^\p)_{k}w(i^\p j^\p)_{k+1}\\
&\qqquad \text{ for some } 1\leq k\leq m(i^\p j^\p), i^\p, j^\p\in T_{1}\\
0,&\text{otherwise}.
}
}
Then $A_{n}=A$ on $T_{1}\times S_{1}$, and $A_{n}(ij)\leq A_{n+1}(ij)\leq A(ij)$ for any $n\geq 1$ and $i,j\in S_{n}$ by the definition. We show that $A_{n}$ is irreducible. Let $i,j\in S_{n}$. We consider the following four cases:\\
Case: $i,j\in T_{n}$. Then $i\cd w(ij)\cd j$ is $A_{n}$-admissible.\\
Case: $i\in T_{n},\ j\in S_{n}\setminus T_{n}$. Then $j=w(i^\p j^\p)_{k}$ for some $k,i^\p,j^\p$. Therefore, $i\cd w(i i^\p)\cd i^\p\cd w(i^\p j^\p)_{1}\cdots w(i^\p j^\p)_{k-1}\cd j$ is $A_{n}$-admissible.\\
Case: $i\in S_{n}\setminus T_{n}$,\ $j\in T_{n}$. Then $i=w(i^\p j^\p)_{k}$ for some $k,i^\p,j^\p$. Thus $i\cd w(i^\p j^\p)_{k+1}\cdots w(i^\p j^\p)_{m(i^\p j^\p)}\cd j^\p\cd w(j^\p j)\cd j$ is $A_{n}$-admissible.\\
Case: $i,j\in S_{n}\setminus T_{n}$. Then $i=w(i^\p j^\p)_{k}$, $j=w(i^\pp j^\pp)_{l}$ for some $k,l,i^\p,j^\p,i^\pp, j^\pp$. We have $i\cd w(i^\p j^\p)_{k+1}\cdots w(i^\p j^\p)_{m(i^\p j^\p)}\cd j^\p\cd w(j^\p i^\pp)\cd i^\pp\cd w(i^\pp j^\pp)_{1}\cdots w(i^\pp j^\pp)_{l-1}\cd j$ is $A_{n}$-admissible. Hence $A_{n}$ is irreducible. 
\proe
We put
$M[S,A]:=\{M\,:\,S_{0}\times S_{0}\to \{0,1\}\,:\,S_{0}\subset S,\ M(ij)\leq A(ij) \text{ for all }i,j\in S_{0}\}$.
The following plays an important role in the proof of our main results.
\thm
{\label{th:P(ph)=supP(phn)}
Let $X$ be a TMS and $\ph\in F^{k+1}(X)$ a real-valued function for a fixed integer $k\geq 1$. Assume that $\ph$ is summable. Then 
\ali
{
P(\ph)=\sup\{P(\ph|_{X_{M}})\,:\,M\in M[S,A] \text{ is irreducible}\}.
}
}
\pros
Since $P(\ph|_{X_{M}})\leq P(\ph)$ for any subsystem $X_{M}$ of $X$, it suffices to show that $P(\ph)\leq \sup_{M\in M[S,A]\,:\,\text{irreducible}}P(\ph|_{X_{M}})$.
If $X$ contains no periodic points, then $P(\ph)=-\infty$ by Theorem \ref{th:pressure_max}. In this case, $P(\ph)=P(\ph|_{X_{M}})=-\infty$ for all irreducible matrices $M\in M[S,A]$.

Now assume that $X$ contains a periodic point of $X$. Then $-\infty<P(\ph)<\infty$. By Theorem \ref{th:pressure_max}, we have $P(\ph)=\max_{T\in S/\leftrightarrow}P(\ph_{X(T)})$ and then $P(\ph)=P(\ph_{X(T)})$ for some $T\in S/\!\!\leftrightarrow$. Therefore we may assume that the transition matrix $A$ of $X$ is irreducible.

Let $S=\{1,2,\dots, d\}$ with $d\leq +\infty$, and let $\hat{A}$ be the matrix indexed by $S\times S$ with all entries equal to $1$. Denote by $\hat{X}$ the full shift with transition matrix $\hat{A}$. Put $N=\{\om\in \hat{X}\,:\,A(\om_{0}\om_{1})=0\}$ and
\ali
{
\ph(\e,\om)=
\case
{
\ph(\om)-(1/\e)\chi_{N}(\om),& \text{if }\ph(\om)\geq \sup_{[\om_{0}]}\ph-1/\e\\
\sup_{[\om_{0}]}\ph-(1/\e)-(1/\e)\chi_{N}(\om),&\text{if }\ph(\om)< \sup_{[\om_{0}]}\ph-1/\e
}
}
for $\e>0$ and $\om\in \hat{X}$. By \cite[Lemma 5.1]{T2024_sum}, conditions (B.1)-(B.3) hold after replacing $A$ with $\hat{A}$ and $M$ with $A$. Moreover, $P(\ph(\e,\cd))\to P(\ph)$ by Theorem \ref{th:perTMS}(3).

On the other hand, consider the finite state space $S_{n}$ and the transition irreducible matrix $A_{n}$ indexed by $S_{n}\times S_{n}$ given in Proposition \ref{prop:inc_irremat}. For convenience, we set $A_{n}(ij)=0$ for $i \in S\setminus S_{n}$ or $j\in S\setminus S_{n}$. 
Define $\hat{M}_{n}\,:\,S\times S\to \{0,1\}$ by $\hat{M}_{n}(ij)=1$ if $i\leq n$ and $\hat{M}_{n}(ij)=0$ if $i\geq n+1$. Let $\ph_{n}(\e,\cd)\,:\,X\to \R$ be
\ali
{
\ph_{n}(\e,\om)=
\case
{
\ph(\e,\om)-1/\e &\text{if }\om_{0}\om_{1}\in S_{n}^{2},\ \ A_{n}(\om_{0}\om_{1})=0,\ A(\om_{0}\om_{1})=1\\
\ph(\e,\om)-1/\e &\text{if }\om_{0}\notin S_{n}\\
\ph(\e,\om) &\text{otherwise}.\\
}
}
Then $\ph_{n}(\e,\cd)$ satisfies also (B.1)-(B.3) replacing $\phe$ by $\ph_{n}(\e,\cd)$ and $d$ by $n$. 
Define $X_{n}$ as the TMS with state space $S_{n}$ and transition matrix $A_{n}$. 
Then $P(\ph_{n}(\e,\cd)_{\hat{M}_{n}})\to P(\ph|_{X_{n}})$ as $\e\to 0$.
Choose any $\eta>0$. By (ii) above, there exists $n_{0}\geq 1$ such that
\ali
{
\sum_{s=n_{0}+1}^{d}\exp(\sup_{\e>0}\sup_{\om\in [s]}\ph(\e,\om))<\eta.
}
Therefore for $n\geq n_{0}$
\ali
{
|(\LR_{\hat{A},\ph(\e,\cd)}f-\LR_{\hat{M}_{n},\ph_{n}(\e,\cd)}f)(\om)|&\leq \sum_{s\in S}|e^{\ph(\e,s\cd \om)}-e^{\ph_{n}(\e,s\cd \om)}|\|f\|_{\infty}\\
&\leq 2\sum_{s=n_{0}+1}^{d}\exp(\ph(s\cd \om))\leq 2\eta\|f\|_{\infty}.
}
Thus $\|\LR_{\hat{A},\ph(\e,\cd)}-\LR_{\hat{M}_{n},\ph_{n}(\e,\cd)}\|_{\infty}\leq 2\eta$
for any $n\geq n_{0}$. Let $(\lam(\e),h(\e,\cd),\nu(\e,\cd))$ be the spectral triplet of $\LR_{\hat{A},\ph(\e,\cd)}$, $(\lam,h,\nu)$ that of $\LR_{A,\ph}$, $(\lam_{n}(\e),h_{n}(\e,\cd),\nu(\e,\cd))$ that of $\LR_{\hat{M}_{n},\ph(\e,\cd)}$, and $(\lam_{n},h_{n},\nu_{n})$ that of $\LR_{A_{n},\ph}$. Note that $\supp\,h_{n}=\bigcup_{s\in S_{n}}[s]$, and $\supp\,\nu_{n}=X_{n}$. Since $A$ is irreducible, $(\lam(\e),h(\e,\cd),\nu(\e,\cd))$ converges to $(\lam,h,\nu)$, and $(\lam_{n}(\e),h_{n}(\e,\cd),\nu_{n}(\e,\cd))$ converges to $(\lam_{n},h_{n},\nu_{n})$. Let $g_{n}(\e,\cd)=h_{n}(\e,\cd)/\|h_{n}(\e,\cd)\|_{\infty}$. Notice the equalities
\ali
{
\lam(\e)=&\exp(P(\phe)),\quad \lam_{n}(\e)=\exp(P(\ph_{n}(\e,\cd)_{\hat{M}_{n}}))\\
\lam=&\exp(P(\ph)),\quad \lam_{n}=\exp(P(\ph_{X_{n}})).
}
Consider the equation
$(\lam(\e)-\lam_{n}(\e))\nu(\e,g_{n}(\e,\cd))=\nu(\e,(\LR_{\hat{A},\phe}-\LR_{\hat{A},\ph_{n}(\e,\cd)})g_{n}(\e,\cd))$ which follows from $\lam(\e)\nu(\e,\cd)=\LR_{\hat{A},\phe}^{*}\nu(\e,\cd)$ and $\lam_{n}(\e)g_{n}(\e,\cd)=\LR_{\hat{A},\ph_{n}(\e,\cd)}g_{n}(\e,\cd)$. Then
\alil
{
\lam(\e)-\lam_{n}(\e)=\frac{\nu(\e,(\LR_{\hat{A},\phe}-\LR_{\hat{A},\ph_{n}(\e,\cd)})g_{n}(\e,\cd))}{\nu(\e,g_{n}(\e,\cd))}.\label{eq:le-lne}
}
Now we show that the left hand side vanishes uniformly in $n$ as $\e\to 0$. To do this, we prove that $\nu(\e,g_{n}(\e,\cd))$ has a uniform positive lower bound in $n$ for sufficiently small $\e>0$. As in the proof of \cite[Proposition 4.6]{T2024_sum}, take a periodic point $\up\in X$ with $\si^{m}\up=\up$ and define $\lam_{*}=(\prod_{i=0}^{m-1}\exp(\ph(\si^{i}\up)))^{1/m}$. Then $\liminf_{\e\to 0}\lam_{n}(\e)\geq \lam_{*}$ for sufficiently large $n$, and so there exists $n_{1}\geq 1$ such that for any $n\geq n_{1}$, $\lam_{n}(\e)^{k}\geq \lam_{*}^{k}/2$. Moreover, Since $\phe$ satisfies $(B.2)$, there exists $n_{2}\geq 1$ such that for any $n>n_{2}$
\ali
{ \sum_{s=n_{2}+1}^{n}\exp(\sup_{[s]}\ph(\e,\cd))<&c_{\adr{us}}^{-k+1}(\lam_{*})^{k}/(6k).
}
Take $\e_{0}>0$ so that for any $0<\e<\e_{0}$
\ali
{
e^{-1/\e}\max_{w\in \{1,\dots, n_{2}\}^{k}}(e^{\sup_{\eta>0}\sup_{[w_{i}]}\ph(\eta,\cd)})^{k}<&\frac{\lam_{*}^{k}}{12 n_{2}^{k}}.
}
We have for $n\geq \max\{n_{0},n_{1},n_{2}\}$
\ali
{
\frac{\lam_{*}^{k}}{2}g_{n}(\e,\om)\leq \lam_{n}(\e)^{k}g_{n}(\e,\om)
= &\LR_{\hat{M}_{n},\ph(\e,\cd)}^{k}g_{n}(\e,\om)\\
\leq&\sum_{w\in \{1,\dots,n_{2}\}^{k}}\exp(S_{k}\ph(\e,w\cd\om))g_{n}(\e,w\cd\om)+\frac{\lam_{*}^{k}}{6}.
}
Take $\tau\in \{1,\dots,n\}^{\Z_{+}}$ so that $g(\e,\tau)>2/3$ and $\tau^{w}\in [w]$ for any $w\in \{1,\dots, n_{2}\}^{k}$. We notice
\ali
{
\exp(S_{k}\ph_{n}(\e,w\cd\tau))
\leq&
\case
{
\exp(S_{k}\ph(\e,w\cd\tau)),&\text{if }w\cd \tau_{0} \text{ is }A_{n}\text{-admissible}\\
\lam_{*}^{k}/(12 n_{2}^{k}),&\text{otherwise}
}
}
Let $c_{\adl{Lc}}=\sup_{\e>0}[\phe]_{k+1}\theta/(1-\theta)$. We see
\ali
{
\frac{\lam_{*}^{k}}{3}<\frac{\lam_{*}^{k}}{2}g(\e,\tau)\leq& \sum_{w\in W_{k}(A_{n})\,:\,w\cd \tau_{0}\in W_{k+1}(A_{n})}\exp(S_{k}\ph(\e,w\cd\tau))g(\e,w\cd\tau)+\frac{\lam_{*}^{k}}{12}+\frac{\lam_{*}^{k}}{6}\\
\leq&c_{\adr{skp}}\sum_{w\in W_{k}(A_{n})\,:\,w\cd \tau_{0}\in W_{k+1}(A_{n})}g(\e,w\cd\tau^{w})+\frac{\lam_{*}^{k}}{4}
}
with $c_{\adl{skp}}=e^{c_{\adr{Lc}}\theta}\max_{w\in \{1,\dots, n_{2}\}^{k}}(e^{\sup_{\eta>0}\sup_{[w_{i}]}\ph(\eta,\cd)})^{k}$, where $W_{k}(A_{n})$ is the set of all $A_{n}$-admissible words with $k$ length.  Since $\supp\,\nu=X$ and $W_{k}(M_{n})\subset W_{k}(A)$, we have $\nu(\e,[w])\to \nu([w])>0$ for any $w\in W_{k}(M_{n})$. Thus, we get
\ali
{
\nu(\e,g_{n}(\e,\cd))\geq& \sum_{w\in \{1,\dots, n_{2}\}^{k}}\nu(\e,\chi_{[w]}g_{n}(\e,\cd))\\
\geq& e^{-c_{\adr{Lc}}\theta}(\min_{v\in \{1,\dots,n_{2}\}^{k}\cap W_{k}(A_{n})}\nu(\e,[v]))\sum_{w\in W_{k}(A_{n})\,:\,w\cd \tau_{0}\in W_{k+1}(A_{n})}g(\e,w\cd\tau^{w})\\
\geq&\frac{\lam_{*}^{k}}{12}e^{-c_{\adr{Lc}}\theta}(\min_{v\in \{1,\dots,n_{2}\}^{k}\cap W_{k}(A_{n})}\nu(\e,[v])).
}
Since the last expression converges to a positive number, there exists $\e_{1}>0$ such that $\nu(\e,g_{n}(\e,\cd))\geq c_{\adl{ngn}}$ for any $0<\e<\min\{\e_{0},\e_{1}\}$ for the positive constant
\ali
{
c_{\adr{ngn}}:=\frac{\lam_{*}^{k}}{24}e^{-c_{\adr{Lc}}\theta}(\min_{v\in \{1,\dots,n_{2}\}^{k}\cap W_{k}(A_{n})}\nu([v])).
}
Consequently, for any $n\geq \max\{n_{0},n_{1},n_{2}\}$ and $0<\e<\min\{\e_{0},\e_{1}\}$, the equation (\ref{eq:le-lne}) implies
\ali
{
|\lam(\e)-\lam_{n}(\e)|\leq \|\LR_{\hat{A},\ph(\e,\cd)}-\LR_{\hat{M}_{n},\ph_{n}(\e,\cd)}\|_{\infty}c_{\adr{ngn}}^{-1}
\leq c_{\adr{ngn}}^{-1}\eta
}
for any $n\geq \max\{n_{0},n_{1},n_{2}\}$ and $0<\e<\min\{\e_{0},\e_{1}\}$. Fix such an $n$. Choose $\e_{2}>0$ so that for $0<\e<\e_{3}$, $|\lam(\e)-\lam|<\eta$ and $|\lam_{n}(\e)-\lam_{n}|<\eta$. Then
\ali
{
|\lam-\lam_{n}|\leq& |\lam-\lam(\e)|+|\lam(\e)-\lam_{n}(\e)|+|\lam_{n}(\e)-\lam_{n}|\leq (2+c_{\adr{ngn}}^{-1})\eta
}
for any $0<\e<\min\{\e_{0},\e_{1},\e_{2}\}$. Hence
\ali
{
\lam=\exp(P(\ph))\leq \lam_{n}+(2+c_{\adr{ngn}}^{-1})\eta\leq \exp(\sup_{M}P(\ph|_{X_{M}}))+(2+c_{\adr{ngn}}^{-1})\eta.
}
Since $\eta>0$ is arbitrary, we conclude $P(\ph)\leq \sup_{M\in M[S,A]\,:\,\text{irreducible}}P(\ph|_{X_{M}})$. 
\proe
\subsection{Formulation of perturbed GIFS with degeneration}\label{sec:GIFSdege}
In this section, we introduce perturbed GIFS in which the images of certain contraction mappings shrink to a single point. Consider the following conditions:
\ite
{
\item[(P.1)] A system $(G=(V,E,i,t),(J_{v}),(O_{v}),(T_{e}))$ is a simple GIFS. Moreover, $G$ has at least one strongly connected component.
\item[(P.2)] For each $\e\in (0,1)$, a system $(\hat{G}=(\hat{V},\hat{E},i,t),(J_{v}(\e))_{v\in \hat{V}},(O_{v}(\e))_{v\in \hat{V}},(T_{e}(\e,\cd))_{e\in \hat{E}})$ is a also simple GIFS . The graph $\hat{G}$ contains $G$ as subgraph. Moreover
\ite
{
\item[(i)] There exists $\e_{0}>0$ such that $c_{\adl{evd}}:=\inf_{0<\e<\e_{0}}\inf_{v\in \hat{V}}\dist (J_{v}(\e),\R^{D}\setminus O_{v}(\e))>0$.
\item[(ii)] There exist $c_{\adl{pcd}}\in (0,c_{\adr{cJO}})$, $\e_{0}>0$ and compact and connected subsets $J_{v}\subset \R^{D}$ ($v\in \hat{V}\setminus V$) with $\sup_{v}\diam J_{v}<+\infty$ such that for any $0<\e<\e_{0}$ and $v\in \hat{V}$, $J_{v}(\e)\subset B(J_{v},c_{\adr{pcd}})$.
\item[(iii)] There exist $a_{v}\in J_{v}$ ($v\in \hat{V}\setminus V$) such that letting $T_{e}(x):\equiv a_{i(e)}$ for $e\in \hat{E}\setminus E$, we have $\sup_{e\in \hat{E}}\sup_{x\in J_{t(e)}(\e)}|T_{e}(\e,x)-T_{e}(x)|\to 0$ as $\e\to 0$.
\item[(iv)] For each $v\in \hat{V}$, $\sup_{e\in \hat{E}\,:\,i(e)=v}\sup_{x\in J_{t(e)}(\e)}|T_{e}^\p(\e,x)-T_{e}^\p(x)|\to 0$ as $\e\to 0$, where $T_{\e}^\p(\e,x)$ means the partial derivative of $T_{\e}(\e,x)$ at $x$.
\item[(vii)] There exists $\e_{0}>0$ such that $r_{1}:=\sup_{0<\e<\e_{0}}\sup_{e\in \hat{E}}\sup_{x\in J_{t(e)}(\e)}|T_{e}^\p(\e,x)|<1$.
\item[(viii)] There exists $c_{\adl{bd3}}>0$ such that $\left|\|T_{e}^\p(\e,x)\|-\|T_{e}^\p(\e,y)\|\right|\leq c_{\adr{bd3}}\|T_{e}^\p(\e,x)\| |x-y|^{\beta}$ for each $e\in \hat{E}$, $x,y\in O_{t(e)}(\e)$.
}
}
We also assume the following as necessary.
\ite
{
\item[(P.3)] There exists $c_{\adl{bd3i}}>0$ such that $\left|\|T_{e}^\p(\e,x)\|_{i}-\|T_{e}^\p(\e,y)\|_{i}\right|\leq c_{\adr{bd3i}}\|T_{e}^\p(\e,x)\|_{i} |x-y|^{\beta}$ for each $e\in \hat{E}$, $x,y\in O_{t(e)}(\e)$.
}
Note that this formulation generalizes the setting considered in \cite{T2019}, which focused on systems defined over finite graphs and assumed $J_{v}(\e)\equiv J_{v}$ and $O_{v}(\e)\equiv O_{v}$. We let
\alil
{
I_{P}^{\hat{V}}:=\{s\geq 0\,:\,\sum_{v\in \hat{V}}\sup_{0<\e<\e_{0}}\sup_{u\in \hat{V}\,:vu\in \hat{E}}\sup_{x\in J_{u}(\e)}\|T_{vu}^\p(\e,x)\|^{s}<+\infty \text{ for some }\e_{0}>0\}.\label{eq:I**=}
}
For convenience, given $\om\in X^{\hat{V}}$, we write $T_{\om_{0}\cdots \om_{n}}$ for the composition $T_{\om_{0}\om_{1}}\circ\cdots \circ T_{\om_{n-1}\om_{n}}$ when $n\geq 1$, and we regard $T_{\om_{0}}$ as the identity map on $O_{\om_{0}}$ when $n=0$. Similarity, we write $T_{\om_{0}\cdots \om_{n}}(\e,\cd)$ for the composition $T_{\om_{0}\om_{1}}(\e,\cd)\circ \cdots\circ T_{\om_{n-1}\om_{n}}(\e,\cd)$ when $n\geq 1$, and define $T_{\om_{0}}(\e,\cd)$ to be the identity map on $O_{\om_{0}}$. Let $\pi(\e,\cd)\,:\,X^{\hat{V}}\to \R^{D}$ be the coding map of the perturbed system. Then the map $\pi\,:\,X^{\hat{V}}\to \R^{D}$ is defined in the same way, namely $\pi\om=\bigcap_{n\geq 1}T_{\om_{0}\cdots \om_{n}}(J_{\om_{n}})$. 
\lem
{\label{lem:conv_pi}
Assume that the conditions (P.1) and (P.2) are satisfied. Then we have $\sup_{\om\in X^{\hat{V}}}|\pi(\e,\om)-\pi(\om)|\to 0$ as $\e\to 0$.
}
\pros
We note $\pi(\e,\si^{k}\om)\in J_{\om_{k}}(\e)\subset B(J_{w_{k}},c_{\adr{pcd}})$ by (P.2)(i). In addition to the fact $c_{\adr{pcd}}<c_{\adr{cJO}}$, we see $B(J_{v},c_{\adr{pcd}}+\eta_{0})\subset O_{v}$ for any $v\in \hat{V}$ for some small $\eta_{0}>0$. For $v\in \ti{V}$, put $U_{v}:=B(J_{v},c_{\adr{pcd}}+\eta_{0})$. 
For $vu\in \hat{E}$, let $\ti{T}_{vu}(\e,\cd):=T_{vu}(\e,\cd)-T_{vu}$. For each $\omega\in X^{\hat{V}}$ and $k\geq 1$, we have
\ali
{
&|\pi(\epsilon,\omega)-\pi(\omega)|=|T_{\omega_{0}\cdots \omega_{k}}(\epsilon,\pi(\epsilon,\sigma^{k}\omega))-T_{\omega_{0}\cdots\omega_{k}}(\pi(\sigma^{k}\omega))|\\
\leq& \sum_{i=0}^{k-1}|T_{\omega_{0}\cdots\omega_{i}}(T_{\omega_{i}\cdots\omega_{k}}(\epsilon,\pi(\epsilon,\sigma^{k}\omega)))-T_{\omega_{0}\cdots\omega_{i+1}}(T_{\omega_{i+1}\cdots\omega_{k}}(\epsilon,\pi(\epsilon,\sigma^{k}\omega)))|\\
&+|T_{\omega_{0}\cdots \omega_{k}}(\pi(\epsilon,\sigma^{k}\omega))-T_{\omega_{0}\cdots\omega_{k}}(\pi(\sigma^{k}\omega))|\\
\leq&\sup_{x\in J_{\om_{1}}(\e)}|\tilde{T}_{\om_{0}\om_{1}}(\e,x)|+c_{\adr{MT}}\sum_{i=1}^{k-1}\sup_{z\in U_{\omega_{i}}}|(T_{\omega_{0}\cdots\omega_{i}})^{\prime}(z)||\tilde{T}_{\omega_{i}\omega_{i+1}}(\epsilon,T_{\omega_{i+1}\cdots\omega_{k}}(\epsilon,\pi(\epsilon,\sigma^{k+1}\omega)))|\\
&\ \ +c_{\adr{MT}}\sup_{z\in U_{\omega_{k}}}|(T_{\omega_{0}\cdots\omega_{k}})^{\prime}(z)||\pi(\epsilon,\sigma^{k}\omega)-\pi(\sigma^{k}\omega)|\\
\leq&\frac{c_{\adr{MT}}}{1-r}\sup_{vu\in \hat{E}}\sup_{x\in J_{u}(\e)}|\tilde{T}_{vu}(\epsilon,x)|+c_{\adr{MT}}\sup_{v\in V}\mathrm{diam}(B(J_{v},c_{\adr{pcd}}))r^{k}
}
from the inequality (\ref{eq:proof_1}) in Proposition \ref{prop:T:mean} by putting $T=T_{\omega_{0}\cdots\omega_{i}}$. Letting $k\to \infty$, $\|\pi(\epsilon,\cdot)-\pi\|_{\infty}\leq {\textstyle \sup_{vu\in \hat{E}}\sup_{x\in J_{u}(\e)}}|\tilde{T}_{vu}(\epsilon,x)|c_{\adr{MT}}/(1-r)$ is satisfied. Thus we obtain the assertion.
\proe
\lem
{\label{lem:lip_pi}
Assume that the conditions (P.1) and (P.2) are satisfied. Then there exists $c_{\adl{cp}}>0$ such that for any $\om,\up\in X^{\hat{V}}$ with $\om_{0}=\up_{0}$, $|\pi(\e,\om)-\pi(\e,\up)|\leq c_{\adr{cp}}d_{r_{1}}(\om,\up)$.
}
\pros
For $\om,\up\in X^{\hat{V}}$ with $\om\neq \up$ and $i(\om_{0})=i(\up_{0})$
, take $k\geq 0$ so that $\om_{0}\cdots\om_{k-1}=\up_{0}\cdots\up_{k-1}$ and $\om_{k}\neq \up_{k}$. Let $\tau=\om_{0}\cdots\om_{k-1}$. We see
\ali
{
|\pi(\e,\om)-\pi(\e,\up)|=&|T_{\tau}(\e,\pi(\e,\si^{k-1}\om))-T_{\tau}(\e,\pi(\e,\si^{k-1}\up))|\\
\leq &c_{\adr{MT}}\sup_{x\in B(J_{\om_{k-1}},c_{\adr{evd}})}|T_{\tau}^\p(\e,x)||\pi(\e,\si^{k-1}\om))|-\pi(\e,\si^{k-1}\up))||\\
\leq &c_{\adr{MT}}(\sup_{v\in \hat{V}}\diam B(J_{v},c_{\adr{pcd}}))r_{1}^{k-2}=c_{\adr{MT}}(\sup_{v\in \hat{V}}\diam B(J_{v},c_{\adr{pcd}}))r_{1}^{-2}d_{r_{1}}(\om,\up).
}
Hence the assertion holds.
\proe
For $\om\in X^{\hat{V}}$, we put
\alil
{
\ph(\e,\om)=\log\|T_{\om_{0}\om_{1}}^\p(\e,\pi(\e,\si\om))\|.\label{eq:phe=...}
}
\lem
{\label{lem:summ}
Assume that the conditions (P.1) and (P.2) are satisfied. Then for any $s\in I_{P}^{\hat{V}}$, $\sum_{v\in \hat{V}}\sup_{0<\e<\e_{0}}\sup_{\om\in [v]}\exp(s\ph(\e,\om))<+\infty$ for some $\e_{0}>0$.
}
\pros
We have the assertion by
$\sup_{\om\in [v]}e^{s\ph(\e,\om)}\leq \sup_{u\,:\,vu\in \hat{E}}\sup_{x\in J_{u}}\|T_{vu}^\p(\e,x)\|^{s}$.
\proe
\lem
{\label{lem:prop_phe}
Assume that conditions (P.1) and (P.2) are satisfied. Then $\sup_{\e>0}[\phe]_{2}<+\infty$ with $\theta=r_{1}^{\beta}$.
}
\pros
For $\om,\up$ with $vu:=\om_{0}\om_{1}=\up_{0}\up_{1}$, 
\ali
{
&|\ph(\e,\om)-\ph(\e,\up)|\\
\leq& \frac{1}{\min(\|T_{vu}^\p(\e,\pi(\e,\si\om))\|,\|T_{vu}^\p(\e,\pi(\e,\si\up)))\|}|\|T_{vu}^\p(\e,\pi(\e,\si\om))\|-\|T_{vu}^\p(\e,\pi(\e,\si\up))\||\\
\leq &c_{\adr{bd3}}|\pi(\e,\si\om)-\pi(\e,\si\up)|^{\beta}\leq c_{\adr{bd3}}r_{1}^{\beta}\theta^{-1}d_{\theta}(\om,\up).
}
Hence we get the assertion.
\proe
For $\om\in X^{\hat{V}}$, we let
$\psi_{s}(\om):=
\case
{
\|T_{\om_{0}\om_{1}}^\p(\pi\si\om)\|^{s},&\text{if } \om_{0}\om_{1}\in E\\
0,&\text{if }\om_{0}\om_{1}\in \hat{E}\setminus E.\\
}$
\lem
{\label{lem:conv_phi_deg}
Assume that the conditions (P.1) and (P.2) are satisfied. Then for each $v\in \hat{V}$ and $s\in I_{P}^{\hat{V}}$, $\sup_{\om\in [v]}|e^{s\ph(\e,\om)}-\psi_{s}(\om)|\to 0$ as $\e\to 0$.
}
\pros
Let $d$ be the minimal integer with $d\geq s$.
Take $v\in \hat{V}$ and $\om\in [v]$. Put $A_{\e}=e^{(s/d)\ph(\e,\om)}$ and $A=\psi_{s}(\om)^{1/D}$. In addition to the facts $A,A_{\e}\leq 1$, we see the inequality $|e^{s\ph(\e,\om)}-\psi_{s}(\om)|=|A_{\e}^{d}-A^{d}|\leq d|A_{\e}-A|$. So it suffices to show that $A_{\e}\to A$ uniformly in $\om\in [v]$. If $\om_{0}\om_{1}\in E$ then
\ali
{
|A_{\e}-A|=&|\|T_{v\om_{1}}^\p(\e,\pi(\e,\si\om))\|^{s/D}-\|T_{v\om_{1}}^\p(\pi\si\om)\|^{s/D}|\\
\leq&|\|T_{v\om_{1}}^\p(\e,\pi(\e,\si\om))\|-\|T_{v\om_{1}}^\p(\pi\si\om)\||^{s/D}\\
\leq&\{|\|T_{v\om_{1}}^\p(\e,\pi(\e,\si\om))\|-\|T_{v\om_{1}}^\p(\pi(\e,\si\om))\| |+|\|T_{v\om_{1}}^\p(\pi(\e,\si\om))\|-\|T_{vu}^\p(\pi\si\om)\||\}^{s/D}\\
\leq&\{\sup_{u\in V\,:\,vu\in E}\sup_{x\in J_{u}(\e)}\|T_{vu}^\p(\e,x)-T_{vu}^\p(x)\| +c_{\adr{bd1}}r\sup_{\up\in X^{\hat{V}}}|\pi(\e,\up)-\pi\up)|^{\beta}\}^{s/D}
}
If $\om_{0}\om_{1}\in \hat{E}\setminus E$ then we have
\ali
{
|A_{\e}-A|=&e^{(s/d)\ph(\e,\om)}=\|T_{v\om_{1}}^\p(\e,x)\|^{s/d}\leq \sup_{u\,:\,vu\in \hat{E}\setminus E}\sup_{x\in J_{u}(\e)}\|T_{vu}^\p(\e,x)\|^{s/d}.
}
Thus $A_{\e}\to A$ uniformly in $\om\in [v]$ by (P.2)(iv).
Hence the proof is valid.
\proe
\thm
{\label{th:conv_phe}
Assume that conditions (P.1) and (P.2) are satisfied. Then for each $s\in I_{P}^{\hat{V}}$, $P(s\ph(\e,\cd))\to P(s\ph)$ as $\e\to 0$.
}
\pros
Together with (\ref{eq:AH=MT}), Theorem \ref{th:perTMS} yields $P(s\phe)=\max_{H\in SC(\hat{G})}P(s\phe|X_{H})$. By condition (P.1), we have $P(s\ph|X_{H})>-\infty$ for some $H\in SC(\hat{G})$. In addition to the assertions of Lemma \ref{lem:summ}, Lemma \ref{lem:prop_phe} and Lemma \ref{lem:conv_phi_deg}, for such a graph $H$, Theorem \ref{th:perTMS}(3) implies $P(s\phe|X_{H})\to P(s\ph|X_{H})=P(s\ph|X_{H}\cap X_{G})$. Since any $H_{0}\in SC(G)$ is a subgraph of $H$ for some $H\in SC(\hat{G})$, we get $P(s\ph)=\max_{H_{0}\in SC(G)}P(s\ph|X_{H_{0}})=\max_{H\in SC(\hat{G})}P(s\ph|X_{H})$. Hence $P(s\phe)$ converges to $P(s\ph)$. 
\proe
Put
\alil
{
&\underline{\ph}(\e,\om):=\log\|T_{\om_{0}\om_{1}}^\p(\e,\pi(\e,\si\om))\|_{i}\quad \text{ for }\om\in X^{\hat{V}}\label{eq:uphe=}\\
&\underline{I}_{P}^{\hat{V}}:=\{s\geq 0\,:\,\sum_{v\in \hat{V}}\sup_{0<\e<\e_{0}}\sup_{u\in \hat{V}\,:vu\in \hat{E}}\sup_{x\in J_{u}(\e)}\|T_{vu}^\p(\e,x)\|_{i}^{s}<+\infty \text{ for some }\e_{0}>0\}.\label{eq:uI**=}
}
Similarity, we obtain the following under additional condition (P.3).
\thm
{\label{th:conv_phei}
Assume that conditions (P.1)-(P.3) are satisfied. Then for each $s\in \underline{I}_{P}^{\hat{V}}$, $P(s\underline{\ph}(\e,\cd))\to P(s\underline{\ph})$ as $\e\to 0$.
}
\section{Main results}\label{sec:dimfor}
\subsection{Upper estimates of the Hausdorff dimension}\label{sec:upper}
We begin with the bounded distortion property. Throughout this section, we consider only the case where $G$ is a simple graph. In fact, any multi GIFS is a special case of a simple GIFS (see Section \ref{sec:reduce_ms}). Let $(G,(J_{v}),(O_{v}),(T_{e}))$ be a simple GIFS. In the following, we set $S:=V$, $A:=A^{V}$, $X:=X^{V}$, $\pi:=\pi^{V}$, and $K:=K^{V}$. 
Let $S^{*}\subset \bigcup_{n=2}^{\infty}V^{n}$ denote the set of all finite paths in $G$ of length at least two. For word $w=w_{1}\cdots w_{k}\in S^{*}$ with $k\geq 2$, let $T_{w}:=T_{w_{1}w_{2}}\circ\cdots \circ T_{w_{k-1}w_{k}}$. For convenience,  define $T_{v}\,:\,O_{v}\to O_{v}$ to be the identity map.
\prop
{\label{prop:BD}
Assume that $(G,(J_{v}),(O_{v}),(T_{e}))$ is a simple GIFS satisfying $(G.7)_{B}$. Let $U_{v}, W_{v}$ be bounded open connected subsets satisfying $U_{v}\subset \overline{U_{v}}\subset W_{v}\subset O_{v}$ for each $v\in V$, such that $\sup_{v}\diam\,U_{v}<\infty$, and $T_{e}U_{t(e)}\subset U_{i(e)}$,  $T_{e}W_{t(e)}\subset W_{i(e)}$ for any $e$. Then:
\ite
{
\item There exists a constant $c_{\adl{BD}}\geq 1$ such that for any finite path $w=w_{1}\cdots w_{k}\in S^{*}$ and any $x,y\in U_{t(w)}$
\ali
{
\textstyle|\log(\prod_{i=1}^{k-1}\|T_{w_{i}w_{i+1}}^\p(T_{w_{i+1}\cdots w_{k}}(x))\|)-\log(\prod_{i=1}^{k-1}\|T_{w_{i}w_{i+1}}^\p(T_{w_{i+1}\cdots w_{k}}(y))\|)|\leq c_{\adr{BD}}|x-y|^{\beta};
}
\item If condition $(G.7)_{BI}$ also holds, then exists a constant $c_{\adl{BD2}}\geq 1$ such that for any finite path $w=w_{1}\cdots w_{k}\in S^{*}$ and any $x,y\in U_{t(w)}$,
\ali
{
\textstyle|\log(\prod_{i=1}^{k-1}\|T_{w_{i}w_{i+1}}^\p(T_{w_{i+1}\cdots w_{k}}(x))^{-1}\|)-\log(\prod_{i=1}^{k-1}\|T_{w_{i}w_{i+1}}^\p(T_{w_{i+1}\cdots w_{k}}(y))^{-1}\|)|\leq c_{\adr{BD2}}|x-y|^{\beta}.
}
}
}
\pros
(1) By condition $(G.7)_{B}$, we have $|T_{e}^\p(y)|\leq (1+c_{\adr{bd1}}|x-y|^{\beta})|T_{e}^\p(x)|$.
Let $x_{i}:=T_{w_{i}\cdots w_{k}}(x)$ and $y_{i}:=T_{w_{i}\cdots w_{k}}(y)$. From Proposition \ref{prop:T:mean} (replacing $J_{v}$ by $\overline{U}_{v}$ and $O_{v}$ by $U_{v}$), we obtain $|x_{i}-y_{i}|\leq c_{\adr{MT}}r^{k-i}|x-y|$. Let $e_{i}:=w_{i}w_{i+1}$. We obtain
\ali
{
&\Big|\prod_{i=1}^{k-1}\|T_{e_{i}}^\p(x_{i+1})\|-\prod_{i=1}^{k-1}\|T_{e_{i}}^\p(y_{i+1})\|\Big|\\
\leq&\sum_{j=1}^{k-1}\prod_{i=1}^{j-1}\|T_{e_{i}}^\p(x_{i+1})\|\left|\|T_{e_{j}}^\p(x_{j+1})\|-\|T_{e_{j}}^\p(y_{j+1})\|\right|\prod_{l=j+1}^{k-1}\|T_{e_{l}}^\p(y_{l+1})\|\\
\leq&\prod_{i=1}^{k-1}\|T_{e_{i}}^\p(x_{i+1})\|\sum_{j=1}^{k-1}c_{\adr{bd1}}|x_{j+1}-y_{j+1}|^\beta\prod_{l=j+1}^{k-1}(1+c_{\adr{bd1}}|x_{l+1}-y_{l+1}|^{\beta})\\
\leq&\prod_{i=1}^{k-1}\|T_{e_{i}}^\p(x_{i+1})\|\sum_{j=1}^{k-1}c_{\adr{bd1}}c_{\adr{MT}}^{\beta}r^{\beta(k-j-1)}|x-y|^\beta\prod_{l=j+1}^{k-1}(1+c_{\adr{bd1}}c_{\adr{MT}}^{\beta}r^{\beta(k-l-1)}|x-y|^{\beta})\\
\leq&c_{\adr{BD}}\prod_{i=1}^{k-1}\|T_{e_{i}}^\p(x_{i+1})\||x-y|^{\beta}
}
by putting $c_{\adr{BD}}=(c_{\adr{bd1}}c_{\adr{MT}}^{\beta}/(1-r^{\beta}))\prod_{i=0}^{\infty}(1+c_{\adr{bd1}}c_{\adr{MT}}^{\beta}\sup_{v}(\diam U_{v})^{\beta})r^{i\beta})$, where the infinite product in $c_{\adr{BD}}$ is convergent by $\sum_{i=0}^{\infty}r^{i\beta}<+\infty$. The desired inequality then follows from the standard estimate $|\log A-\log B|\leq |A-B|/\min(A,B)$ with $A:=\prod_{i=1}^{k-1}\|T_{e_{i}}^\p(x_{i+1})\|$ and $B:=\prod_{i=1}^{k-1}\|T_{e_{i}}^\p(y_{i+1})\|$.
\smallskip
\\
(2) By combining Proposition \ref{prop:T:mean_inf} with a similar argument as in (1), we obtain the desired result
with the constant $c_{\adr{BD2}}=(c_{\adr{bd2}}c_{\adr{iMT}}^{\beta}/(1-r^{\beta}))\prod_{i=0}^{\infty}(1+c_{\adr{bd2}}c_{\adr{iMT}}^{\beta}\sup_{v}(\diam U_{v})^{\beta})r^{i\beta})$.
\proe
\cor
{\label{cor:BD_v2}
Assume that a set $(G,(J_{v}),(O_{v}),(T_{e}))$ is a simple GIFS and that condition $(G.7)_{B}$ is satisfied. Let $U_{v}$ and $W_{v}$ be bounded, open, connected subsets such that $U_{v}\subset \overline{U_{v}}\subset W_{v}\subset O_{v}$ $(v\in V)$ with $\sup_{v}\diam U_{v}<\infty$, and for any $e\in E$, $T_{e}U_{t(e)}\subset U_{i(e)}$,\ $T_{e}W_{t(e)}\subset W_{i(e)}$. Then:
\ite
{
\item There exists $c_{\adl{BDv}}\geq 1$ such that for any finite path $w=w_{1}\cdots w_{k}\in S^{*}$ and any $x,y\in U_{t(w)}$, $c_{\adr{BDv}}^{-1}\leq \prod_{i=1}^{k-1}\|T_{w_{i}w_{i+1}}^\p(T_{w_{i+1}\cdots w_{k}}(x))\|/\prod_{i=1}^{k-1}\|T_{w_{i}w_{i+1}}^\p(T_{w_{i+1}\cdots w_{k}}(y))\|\leq c_{\adr{BDv}}$.
\item If $(G.7)_{BI}$ also holds, then there exists $c_{\adl{BDv2}}\geq 1$ such that for any $w\in S^{*}$ and any $x,y\in U_{t(w)}$, $c_{\adr{BDv2}}^{-1}\leq \prod_{i=1}^{k-1}\|T_{w_{i}w_{i+1}}^\p(T_{w_{i+1}\cdots w_{k}}(x))\|_{i}/\prod_{i=1}^{k-1}\|T_{w_{i}w_{i+1}}^\p(T_{w_{i+1}\cdots w_{k}}(y))\|_{i}\leq c_{\adr{BDv2}}$.
}
}
\pros
The assertions follow from Proposition \ref{prop:BD} by letting $c_{\adr{BDv}}=\exp(c_{\adr{BD}}\sup_{v}(\diam U_{v})^{\beta})$ and $c_{\adr{BDv2}}=\exp(c_{\adr{BD2}}\sup_{v}(\diam U_{v})^{\beta})$.
\proe
\thm
{\label{th:GIFS_upper}
Assume that $(G,(J_{v}),(O_{v}),(T_{e}))$ is a simple GIFS, that $G$ is strongly connected, and that conditions $(G.7)_{B}$ and $(G.9)_{J}$ are satisfied.
Then $\dim_{H}K\leq \overline{s}$, where $\overline{s}:=\inf\{s\geq 0\,:\,P(s\ph^{V})\leq 0\}$.
}
\pros
Let $\ph:=\ph^{V}$ for simplicity. Choose small $\delta>0$ such that $\delta<c_{\adr{cJO}}$, and define $U_{v}:=\bigcup_{z\in J_{v}}B(z,\delta)$,\  $W_{v}:=\bigcup_{z\in J_{v}}B(z,\delta_{1})$ for some $\delta<\delta_{1}<c_{\adr{cJO}}$. Since each $T_{e}$ is a contraction on $O_{t(e)}$, it follows that $T_{e}U_{t(e)}\subset U_{i(e)}$ and $T_{e}W_{t(e)}\subset W_{i(e)}$. Fix any $t>\inf\{s\geq 0\,:\,P(s\ph)\leq 0\}$. For each $w\in S^{n}$, choose $z(w)\in \hat{J}_{w_{n}}$ and for each $v\in V$, choose $y(v)\in J_{v}$. Then we have
\ali
{
0>&P(t\ph)=\lim_{n\to \infty}\frac{1}{n}\log\sum_{w\in S^{n}}\exp(\sup_{\om\in [w]}S_{n}(t\ph(\om)))\\
=&\lim_{n\to \infty}\frac{1}{n}\log\sum_{w\in S^{n}}\sup_{\om\in [w]}\Big(\prod_{k=1}^{n-1}\|T_{w_{k}w_{k+1}}^\p(\pi\si^{k}\om)\|^{t}\Big)\|T_{w_{n}\om_{n}}^\p(\pi\si^{n}\om)\|^{t}\\
=&\lim_{n\to \infty}\frac{1}{n}\log\sum_{w\in S^{n}}\sup_{\om\in [w]}\Big(\prod_{k=1}^{n-1}\|T_{w_{k}w_{k+1}}^\p(T_{w_{k+1}\cdots w_{n}}(\pi\si^{n-1}\om))\|^{t}\Big)\|T_{w_{n}\om_{n}}^\p(\pi\si^{n}\om)\|^{t}\\
\geq&\limsup_{n\to\infty}\frac{1}{n}\log\Big((c_{\adr{BDv}})^{2}\sum_{w\in S^{n}}\Big(\prod_{k=1}^{n-1}\|T_{w_{k}w_{k+1}}^\p(T_{w_{k+1}\cdots w_{n}}(y(w)))\|^{t}\Big)\sup_{v\in V\,:\,\atop{w_{n}v\in E}}\|T_{w_{n}v}^\p(z(v))\|^{t}\Big)
}
by using Corollary \ref{cor:BD_v2}(1). Therefore there exists $n_{0}\geq 2$ such that for any $n\geq n_{0}$, the following holds:
\alil
{
\sum_{w\in S^{n}}\Big(\prod_{k=1}^{n-1}\|T_{w_{k}w_{k+1}}^\p(T_{w_{k+1}\cdots w_{n}}(y(w)))\|^{t}\Big)\sup_{v\in V\,:\,\atop{w_{n}v\in E}}\|T_{w_{n}v}^\p(z(v))\|^{t}\leq \exp(nP(t\ph)/2).\label{eq:sumleq exp(nP(tph)/2)}
}

Now take $\delta>0$ as above, and choose $n\geq n_{0}$ such that $c_{\adr{MT}}r^{n}<\delta$. For any $A$-admissible word $w\in S^{n}$, by Proposition \ref{prop:T:mean}, we have $\diam T_{w}J_{t(w)}\leq c_{\adr{MT}}\sup_{z\in U_{t(w)}}\|T_{w}^\p(z)\|\leq c_{\adr{MT}}r^{n}<\delta$ is satisfied. Therefore, the collection $(T_{w}J_{t(w)})_{w\in S^{n}}$ forms a closed $\delta$-cover of $K$. By using bounded distortion property, we have
\ali
{
\mathcal{H}^{t}_{\delta}(K):=&\inf\{\sum_{k=1}^{\infty}(\diam C_{k})^{t}\,:\,\{C_{k}\} \text{ is a closed cover of } K \text{ and }\diam C_{k}< \delta\}\\
\leq& \sum_{w\in S^{n}}\diam(T_{w}J_{t(w)})^{t}\\
\leq &c_{\adr{MT}}^{t}\sum_{w\in S^{n}}\sup_{z\in U_{t(w)}}\|T_{w}^\p(z)\|^{t}\diam(J_{t(w)})^{t}\\
\leq&c_{\adr{MT}}^{t}\sum_{w\in S^{n}}\sup_{z\in U_{t(w)}}\prod_{i=1}^{n}\|T_{w_{i}w_{i+1}}^\p(T_{w_{i+1}\dots w_{n}}(z))\|^{t}\diam(J_{w_{n}})^{t}\\
\leq&c_{\adr{CJ}}^{t}c_{\adr{MT}}^{t}c_{\adr{BDv}}^{2t}\sum_{w\in S^{n}}\prod_{i=1}^{n}\|T_{w_{i}w_{i+1}}^\p(T_{w_{i+1}\dots w_{n}}(y(w)))\|^{t}\sup_{v\in V\,:\,w_{n}v\in E}\sup_{z\in J_{v}}\|T_{w_{n}v}^\p(z(v))\|^{t}\\
\leq&c_{\adr{CJ}}^{t}c_{\adr{MT}}^{t}c_{\adr{BDv}}^{2t}\exp(nP(t\ph)/2)\qquad (\because (\ref{eq:sumleq exp(nP(tph)/2)}))
}
for any $n\geq n_{0}$. Letting $n\to \infty$, we observe that $\exp(nP(t\ph)/2)\to 0$, and hence $\mathcal{H}^{t}_{\delta}(K)=0$ for any $\delta>0$. Thus $\mathcal{H}^{t}(K)=\lim_{\delta\to 0}\mathcal{H}^{t}_{\delta}(K)=0$, which implies $\dim_{H}K\leq t$. Since $t>\overline{s}$ is arbitrary, we conclude that $\dim_{H}K\leq \overline{s}$.
\proe
Finally, we provide an upper estimate in the general case where $G$ is an arbitrary graph and condition $(G.8)_{BV}$ is satisfied. We define the constant $c_{\adl{scV2}}$ by
\alil
{
c_{\adr{scV2}}:=&\inf\{s\geq 0\,:\,\sum_{v\in V}\sup_{v^\p\in V\,:\,vv^\p\in E}\sup_{x\in K\cap J_{v^\p}}\|T_{vv^\p}(x)^\p\|^{s}<\infty\}.
}
\thm
{\label{th:GIFS_upper_anygraph}
Assume that $(G,(J_{v}),(O_{v}),(T_{e}))$ is a simple GIFS and that conditions $(G.7)_{B}$ and $(G.9)_{J}$ are satisfied. Then $\dim_{H}K\leq \max(\overline{s},c_{\adr{scV2}})$, where $\overline{s}:=\inf\{s\geq 0\,:\,P(s\ph^{V})\leq 0\}$.
}
\pros
If $\dim_{H}K=0$, then the inequality clearly holds. Therefore, we may assume that $\dim_{H}K>0$.
To prove the result, we construct a perturbed simple GIFS $\mathcal{G}(\e):=(\hat{G},(J_{v}(\e)), (O_{v}(\e)),(T_{e}(\e,\cd)))$ satisfying (P.1) and (P.2) in Section \ref{sec:GIFSdege}:
\ite
{
\item Let $\hat{V}:=V$, $\hat{E}:=V^{2}$, $J_{v}(\e):=J_{v}$ and $O_{v}(\e):=O_{v}$. In particular, $\hat{G}=(\hat{V},\hat{E})$ is strongly connected.
\item For $e\in E$, $T_{e}(\e,\cd):=T_{e}$.
\item For $e=vu\in \hat{E}\setminus E$, we take $x_{vu}\in \mathrm{int} J_{v}$ and $T_{vu}(\e,\cd)$ so that
\ite
{
\item[(i)] $T_{vu}(\e,x)$ has a conformal mapping $T_{vu}(\e,x):=r_{vu}(\e)A_{vu}(x)+b_{vu}(\e)$ with $0<r_{vu}(\e)<1$, $A_{e}$ is a linear isometry on $\R^{D}$ and $b_{vu}(\e)\in \R^{D}$.
\item[(ii)] $T_{vu}(\e, J_{u})\subset B(x_{vu},\eta)\subset J_{v}$ and $\sup_{vu}\sup_{x\in J_{u}}|T_{vu}(\e,x)-x_{vu}|\to 0$ as $\e\to 0$.
\item[(iii)] $\sum_{vu\in \hat{E}\setminus E}\sup_{\e>0}r_{vu}(\e)^{t}<+\infty$ for any $t>0$ and $\sup_{vu\in \hat{E}\setminus E}r_{vu}(\e)\to 0$ as $\e\to 0$.
}
}
Note that even if $T_{e}(\e,O_{t(e)})\not\subset O_{i(e)}$, we may redefine $O_{v}:=B(J_{v},c_{\adr{cJO}}/2)$ for each $v$, so that $T_{e}(\e,O_{t(e)})\subset O_{i(e)}$ for all $e\in \hat{E}$. It is then easy to verify that $\mathcal{G}(\e)$ satisfies conditions $(P.1)$ and $(P.2)$. Moreover, by (2)(iii) above, any $s>c_{\adr{scV2}}$ belongs to the interval $I_{P}^{\hat{V}}$ defined in (\ref{eq:I**=}).

Now choose any $s>\max(\overline{s},c_{\adr{scV2}})$. By Theorem \ref{th:conv_phe}, we have $P(s\phe)\to P(s\ph^{V})<0$ as $\e\to 0$. Let $K(\e)$ denote the limit set of the perturbed GIFS $\mathcal{G}(\e)$. Then by Theorem \ref{th:GIFS_upper}, we see $\dim_{H}K(\e)\leq s$. Since $K\subset K(\e)$, it follows that $\dim_{H}K\leq s$. As $s>\max(\overline{s},c_{\adr{scV2}})$ was arbitrary, the conclusion follows.
\proe
\subsection{Lower estimates of the Hausdorff dimension}\label{sec:lower}
In this section, we use the notation $S,A,X,\pi,K,S^{*}$ as defined in Section \ref{sec:upper}.
\prop
{\label{prop:prop_sepa}
Assume that $(G,(J_{v}),(O_{v}),(T_{e}))$ is a simple GIFS, that $G$ is strongly connected, and that conditions $(G.6)_{S}$ and $(G.7)_{BI}$ are satisfied. Let $W_{v}\subset J_{v}$ be non-empty subsets for each $v\in V$ such that $T_{e}W_{t(e)}\subset W_{i(e)}$ for all $e\in E$. Define $\Delta(v,u,u^\p):=\dist(T_{vu}W_{u},T_{vu^\p}W_{u^\p})$ for $vu,vu^\p\in E$ with $u\neq u^\p$. Then there exists a constant $c_{\adl{sep}}>0$ such that for any finite path $w=w_{1}\cdots w_{k}\in S^{*}$, any $u,u^\p\in E$ with $w_{k}u,w_{k}u^\p\in E$, and any $z\in J_{w_{k}}$, we have $\dist (T_{wu}W_{u},T_{wu^\p}W_{u^\p})\geq c_{\adr{sep}}\Delta(w_{k},u,u^\p)\prod_{i=1}^{k}\|T_{w_{i}w_{i+1}}^\p(T_{w_{i+1}\cdots w_{k}}(z))\|_{i}$.
}
\pros
For each $v\in V$, let $U_{v}:=B(J_{v},\delta)\subset O_{v}$ for sufficiently small $\delta>0$. Then each $U_{v}$ is bounded, open, and connected, and satisfies $T_{e}U_{t(e)}\subset U_{i(e)}$ and $W_{v}\subset U_{v}$. We have
\ali
{
\dist (T_{wu}W_{u},T_{wu^\p}W_{u^\p})=&\inf_{x\in W_{u},\ y\in W_{u^\p}}|T_{w}(T_{w_{k}u}(x))-T_{w}(T_{w_{k}u^\p}(y))|\\
\geq&c_{\adr{iMT}}\inf_{a\in U_{w_{k}}}\|T_{w}^\p(a)\|_{i}\inf_{x\in W_{u},\ y\in W_{u^\p}}|T_{w_{k}u}(x)-T_{w_{k}u^\p}(y)|\\
\geq&c_{\adr{iMT}}\inf_{a\in U_{t(e)}}\prod_{i=1}^{k}\|T_{w_{i}w_{i+1}}^\p(T_{w_{i+1}\cdots w_{k}}(a))\|_{i}\Delta(e,e^\p)\\
\geq&c_{\adr{iMT}}c_{\adr{BDv2}}^{-1}\prod_{i=1}^{k}\|T_{w_{i}w_{i+1}}^\p(T_{w_{i+1}\cdots w_{k}}(z))\|_{i}\Delta(w_{k},u,u^\p)
}
by Corollary \ref{cor:BD_v2} and Proposition \ref{prop:T:mean_inf}.
Hence, the assertion follows with $c_{\adr{sep}}:=c_{\adr{iMT}}c_{\adr{BDv2}}^{-1}$.
\proe
For each $v\in V$, define $K_{v}:=\{\om\in K\,:\,i(\om)=v\}$. Then $K=\bigcup_{v\in V}K_{v}$, $T_{e}K_{t(e)}\subset K_{t(e)}$, and $K_{v}\subset J_{v}$.
We next state a proposition that requires the graph to be finite.
\prop
{\label{prop:GIFS_lower_hole_ball}
Assume that $(G,(J_{v}),(O_{v}),(T_{e}))$ is a finite and simple graph GIFS, $G$ is strongly connected, and that conditions $(G.6)_{S}$ and $(G.7)_{BI}$ are satisfied. Then there exists a constant $c_{\adl{lowhb}}>0$, for any $\kappa\in (0,1)$, $k\geq 2$, any path $w=w_{1}\cdots w_{k}\in S^{*}$, any $z\in J_{w_{k}}$, and any $x\in T_{w}K_{w_{k}}$, we have $K\cap B(x,c_{\adr{lowhb}}\kappa \prod_{i=1}^{k}\|T_{w_{i}w_{i+1}}^\p(T_{w_{i+1}\cdots w_{k}}(z))\|_{i})\subset T_{w}K_{w_{k}}$.
}
\pros
Choose any distinct paths $w,\tau\in S^{k}$. Let $1\leq l\leq k$ be such that $\gamma:=w_{1}\cdots w_{l-1}=\tau_{1}\cdots \tau_{l-1}$ and $w_{l}\neq \tau_{l}$. Applying Proposition \ref{prop:prop_sepa} with $W_{v}$ replaced by $K_{v}$, we obtain
\ali
{
\dist(T_{w}K_{w_{k}}, T_{\tau}K_{\tau_{k}})\geq&\dist(T_{\gamma\cd w_{l}}K_{w_{l}},T_{\gamma\cd \tau_{l}}K_{\tau_{l}})\\
\geq&c_{\adr{lowhb}}\prod_{i=1}^{l-1}\|T_{w_{i}w_{i+1}}^\p(T_{w_{i+1}\cdots w_{l-1}w_{l}}(T_{w_{l}\cdots w_{k}}(z)))\|_{i}\\
\geq&c_{\adr{lowhb}}\prod_{i=1}^{k}\|T_{w_{i}w_{i+1}}^\p(T_{w_{i+1}\cdots w_{k}}(z))\|_{i}\geq c_{\adr{lowhb}}\kappa\prod_{i=1}^{k}\|T_{w_{i}w_{i+1}}^\p(T_{w_{i+1}\cdots w_{k}}(z))\|_{i}
}
by putting $c_{\adr{lowhb}}=c_{\adr{sep}}\min\{\Delta(v,u,u^\p)\,:\,v,u,u^\p\in V,\ vu,vu^\p\in E,\ u\neq u^\p\}$. Due to condition $(G.6)_{S}$ and the finiteness of $S=V$, we have $\Delta(v,u,u^\p)>0$, and therefore $c_{\adr{lowhb}}>0$. Thus
\ali
{
T_{\gamma}K_{\gamma_{k}} \cap B(x,c_{\adr{lowhb}}\kappa\prod_{i=1}^{k}\|T_{w_{i}w_{i+1}}^\p(T_{w_{i+1}\cdots w_{k}}(z))\|_{i})=\emptyset
}
for any $x\in T_{w}W_{w_{k}}$. Hence the assertion follows from the fact $K=\bigcup_{w\in E^{k}}T_{w}K_{t(w)}$.
\proe
\prop
{\label{prop:low_f->inf}
Assume that $\mathcal{G}:=(G,(J_{v}),(O_{v}),(T_{e}))$ is a simple GIFS satisfying condition $(G.8)_{SV}$. Also assume that for any finite strongly connected subgraph $G(0):=(V(0),E(0))$ of $G$, the subsystem $(G(0)=(V(0),E(0)),(J_{v})_{v\in V(0)},(O_{v})_{v\in V(0)},(T_{e})_{e\in E(0)})$ admits a lower estimate: taht is, its limit set $K(0)$ satisfies $\dim_{H}K(0)\geq \underline{s}(0)$, where $P(\underline{s}(0)\underline{\ph}_{0})=0$ and $\underline{\ph}_{0}(\om):=\log\|T_{\om_{0}\om_{1}}^\p(\pi\si\om)\|_{i}$ for $\om\in X_{G(0)}$. Then the inequality $\underline{s}\leq \dim_{H}K$ holds, where $\underline{s}:=\inf\{s\geq 0\,:\,P(s\underline{\ph})\leq 0\}$.
}
\pros
We divide the proof into two cases:
\smallskip
\\
\textbf{Case I}: $G$ is strongly connected. Choose a sequence of subgraphs $G_{n}=(V_{n},E_{n})$ of $G$ $(n\geq 1)$ such that each $G_{n}$ is finite and strongly connected, $V_{n}\subset V_{n+1}$, $E_{n}\subset E_{n+1}$, and $\bigcup_{n}V_{n}=V$,\ $\bigcup_{n}E_{n}=E$. (see \cite[Lemma 2.7.2]{MU} for the existence of such a sequence.)
Define $\underline{\ph}_{n}(\om):=\log\|T_{\om_{0}\om_{1}}^\p(\pi\si\om)\|_{i}$ for $\om\in X_{G_{n}}$, and let $\alpha_{n}>0$ be such that $P(\alpha_{n}\underline{\ph}_{n})=0$.
By assumption, the limit set $K^{n}$ of the subsystem $(G_{n},(J_{v})_{v\in V_{n}},(O_{v})_{v\in V_{n}}, (T_{e})_{e\in E_{n}})$ satisfies $\dim_{H}K^{n}\geq \alpha_{n}$. From the monotonicity of pressure, we have
$0=P(\alpha_{n}\underline{\ph}_{n})\leq P(\alpha_{n}\underline{\ph}_{n+1})$
and therefore $\alpha_{n}\leq \alpha_{n+1}$. Let $\alpha_{\infty}:=\lim_{n\to \infty}\alpha_{n}$. Observe the inequality
\alil
{
\alpha_{\infty}\leq \lim_{n\to \infty}\dim_{H}K^{n}\leq \dim_{H}K.\label{eq:ainf<=dimK}
}
Now we verify that $\underline{s}\leq \alpha_{\infty}$. Fix any $\delta>0$ and $n\geq 1$. We see
\ali
{
P((\alpha_{\infty}+\delta)\underline{\ph}_{n})=&\lim_{m\to \infty}\frac{1}{m}\log\sum_{w\in V_{n}^{m}}\sup_{\om\in [w]}\exp(((\alpha_{\infty}+\delta)S_{m}\underline{\ph}_{n}(\om)))\\
\leq&\delta\log r+P(\alpha_{\infty}\underline{\ph}_{n})\leq \delta\log r<0.
}
By condition $(G.8)_{SV}$ and the continuity of $s\mapsto P(s\underline{\ph}_{n})$, there exists $\delta_{0}>0$ such that for any $0<\delta_{1}<\delta_{0}$, we have $P((\underline{s}+\delta_{1})\underline{\ph}_{n})>\delta\log r$. Since $(\underline{s}+\delta_{1})\underline{\ph}$ is summable, Theorem \ref{th:P(ph)=supP(phn)} implies that
$\lim_{n\to \infty}P((\underline{s}+\delta_{1})\underline{\ph}_{n})=P((\underline{s}+\delta_{1})\underline{\ph})>\delta\log r$.
Thus there exists $n\geq 1$ such that
\ali
{
P((\underline{s}+\delta_{1})\underline{\ph}_{n})>\delta\log r\geq P((\alpha_{\infty}+\delta)\underline{\ph}_{n}).
}
This yields the inequality $\underline{s}<\underline{s}+\delta_{1}<\alpha_{\infty}+\delta$ for any $\delta>0$ and so we get $\underline{s}\leq \alpha_{\infty}$. 
Combining this with the inequality $\alpha_{\infty}\leq \dim_{H} K$, we conclude that $\underline{s}\leq \dim_{H} K$.
\smallskip
\\
\textbf{Case II}: $G$ is an arbitrary graph. By condition $(G.8)_{SV}$, the function $\underline{s}\underline{\ph}^{V}$ is summable. Applying Theorem \ref{th:pressure_max}, we obtain $0=P(\underline{s}\underline{\ph}^{V})=\max_{T\in S/\leftrightarrow}P(\underline{s}\underline{\ph}^{V}|X_{M(T)})$. Therefore, there exists a strongly connected component $G_{0}=(V_{0},E_{0})$ of $G$ such that $0=P(\underline{s}\underline{\ph}^{V_{0}})$. By Case I, we then have $\underline{s}\leq \dim_{H}K^{\hat{V}}\leq \dim_{H}K$. Hence the proof is complete.
\proe
Now we are in a position to state a lower estimate for the Hausdorff dimension.
\thm
{\label{th:GIFS_lower}
Assume that $(G,(J_{v}),(O_{v}),(T_{e}))$ is a simple GIFS satisfying conditions $(G.6)_{S}$, $(G.7)_{B}$, $(G.7)_{BI}$, and $(G.8)_{SV}$. Then $\underline{s}\leq \dim_{H}K$, where $\underline{s}:=\inf\{s\geq 0\,:\,P(s\underline{\ph}^{V})\leq 0\}$.
}
\pros
FIrst, suppose that $G$ contains no strongly connected components. Then, for any $T\in S/\!\!\leftrightarrow$, the matrix $M(T)$ is a $1\times 1$ zero matrix by equation (\ref{eq:AH=MT}). In this case, condition (B.2) fails, and thus Theorem \ref{th:pressure_max} gives $P(s\ph)=-\infty$ for any $s>0$. Thus $\underline{s}=0$, and the desired inequality $\underline{s}\leq \dim_{H} K$ is trivially satisfied.
Now assume that $G$ has at least one strongly connected component.
We divide the proof into two cases.
\smallskip
\\
\textbf{Case I}: $G$ is finite strongly connected graph. Choose $\alpha>0$ such that $P(\alpha \underline{\ph}^{V})=0$. Let $\mu$ be the Gibbs measure for the potential $\alpha\underline{\ph}^{V}$, and define $\ti{\mu}=\mu\circ \pi^{-1}$. Let $A\subset K$ be a non-empty Borel set, and fix a point $x\in A$ with $\om=\pi^{-1}x$. We consider the ratio $\ti{\mu}(B(x,r))/r^{\alpha}$ for small $r>0$. Let $\kappa:=(\min_{e\in E}\inf_{z\in J_{t(e)}}\|T_{e}^\p(z)\|_{i})^{-1}$, and choose $0<r<\kappa^{-1}$. For each $n\geq 0$, define $r_{n}=\prod_{i=0}^{n}\|T_{\om_{i}\om_{i+1}}^\p(T_{\om_{i+1}\cdots \om_{n}}(\pi\si^{n+1}\om))\|_{i}$. Then there exists a unique integer $k=k(r)>0$ such that $r_{k}<r\leq r_{k-1}$, since $r_{k-1}\leq r_{k}\kappa^{-1}$. Therefore we have $\kappa r_{k}<\kappa r\leq r_{k}$. By virtue of Proposition \ref{prop:GIFS_lower_hole_ball},
\ali
{
&K\cap B(x,c_{\adr{lowhb}}\kappa \prod_{i=0}^{k}\|T_{\om_{i}\om_{i+1}}^\p(T_{\om_{i+1}\cdots \om_{k}}(\pi\si^{k+1}\om))\|_{i})\subset T_{w}K_{t(w)}
\subset T_{\om_{0}\cdots \om_{k}}(K_{t(\om_{k})}),
}
where each set $J_{v}$ is given in Proposition \ref{prop:GIFS_lower_hole_ball}. Thus 
\ali
{
&\ti{\mu}(B(x, c_{\adr{lowhb}}\kappa r))\leq \ti{\mu}(B(x,c_{\adr{lowhb}}r_{k}))\leq\ti{\mu}(T_{\om_{0}\cdots \om_{k}}(K_{t(\om_{k})}))
=\mu([\om_{0}\cdots \om_{k}])\\
\leq &c_{\adr{Gibbs1}}e^{\alpha S_{k+1}\underline{\ph}(\om)}
= c_{\adr{Gibbs1}}(\prod_{i=0}^{k}\|T_{\om_{i}}^\p(\pi\si^{i+1}\om)\|_{i})^{\alpha}
= c_{\adr{Gibbs1}} (r_{k})^{\alpha}<c_{\adr{Gibbs1}}r^{\alpha}=c_{\adr{Gibbs1}}(c_{\adr{lowhb}}\kappa)^{-\alpha}(c_{\adr{lowhb}}\kappa r)^{\alpha},
}
where $c_{\adr{Gibbs1}}$ is a constant appearing in the definition of Gibbs measure (\ref{eq:Gibbs}).
By Frostman lemma, the $\alpha$-Hausdorff measure $\mathcal{H}^{\alpha}$ satisfies
\ali
{
\mathcal{H}^{\alpha}(K)\geq \ti{\mu}(K)/(c_{\adr{Gibbs1}}(c_{\adr{lowhb}}\kappa)^{-\alpha})>0.
}
This means that $\alpha\leq \dim_{H}K$. We obtain the assertion under a finite graph case.
\smallskip
\\
\textbf{Case II}: $G$ is an arbitrary graph. By Proposition \ref{prop:low_f->inf} and Case I, we get the assertion. Hence the proof is complete.
\proe
\subsection{Bowen's formula}\label{sec:Bowen}
We now state and prove a Bowen equation for conformal mappings associated with graph-directed systems.
\thm
{\label{th:GIFS_Bowen}
Assume that $(G,(J_{v}),(O_{v}),(T_{e}))$ is a simple GIFS, and that conditions $(G.5)_{C}$, $(G.7)_{B}$, $(G.8)_{SV}$ and $(G.9)_{J}$ are satisfied. Assume also that either the SSC $(G.6)_{S}$ or the OSC $(G.6)_{O}$ holds. Then $\dim_{H}K=s$ if and only if $P(s\ph)=0$, where $\ph:=\ph^{V}=\underline{\ph}^{V}$. In particular, $\dim_{H}K=\max_{G_{0}}\dim_{H}K(G_{0})$, where the maximum is taken over all strongly connected components $G_{0}$ of $G$, and $K(G_{0})$ denote the limit set of the subsystem $(G_{0},(J_{v}),(O_{v}),(T_{e}))$.
}
\pros
Since $\|T_{e}^\p(x)\|=\|T_{e}^\p(x)\|_{i}$ and $(G.8)_{SV}$ is satisfied, we have  $c_{\adr{scV2}}\leq \overline{s}$. Therefore the upper bound $\dim_{H}\leq s$ follows from Theorem \ref{th:GIFS_upper}. 
To prove the lower bound, first assume that the SSC $(G.6)_{S}$ holds. Then $\dim_{H}K\geq s$ is guaranteed by Theorem \ref{th:GIFS_lower}.

Next, suppose that the OSC $(G.6)_{O}$ holds. Assume also that $G$ is finite and strongly connected graph. In this case, the general dimension theory for conformal mapping under the OSC yields $\dim_{H}K=s$ (see, e.g., \cite[Theorem 19.6.3.]{URM}). Proposition \ref{prop:low_f->inf} extends this lower estimate, completing the proof of the equality.

Finally, we prove the last assertion. Clearly, for any strongly connected component $G_{0}$, we have $\dim_{H}K(G_{0})\leq \dim_{H}K$. Since $s\ph$ is summable, we obtain $0=P(s\ph)=\max_{T\in V/\leftrightarrow}P(s\ph|X_{M(T)})=\max_{G_{0}}P(s\ph|X^{V_{0}})$. Thus, there exists a component $G_{0}$ such that $0=P(s\ph|X^{V_{0}})$, implying $\dim_{H} K(G_{0})=s$. Hence $\dim_{H} K=\max_{G_{0}}\dim_{H}K(G_{0})$.
\proe
As a direct consequence, we obtain Bowen formula for multi GIFS.
\cor
{\label{cor:GIFS_Bowen}
Assume that $(G,(J_{v}),(O_{v}),(T_{e}))$ is a multi GIFS, and that conditions $(G.5)_{C}$, $(G.7)_{B}$ and $(G.8)_{SE}$ are satisfied. Assume also that either $(G.6)_{S}$ or $(G.6)_{O}$ holds. Then $\dim_{H}K=s$ if and only if $P(s\ph)=0$, where $\ph:=\ph^{E}=\underline{\ph}^{E}$. In particular, $\dim_{H}K=\max_{G_{0}}\dim_{H}K(G_{0})$.
}
\pros
This follows directly from Proposition \ref{prop:reduce_GIFS} and Theorem \ref{th:GIFS_Bowen}.
\proe
\subsection{Continuity of the dimension in perturbed nonconformal GIFS with degeneration}\label{sec:conv_dege}
Recall conditions (P.1), (P.2) and (P.3) from Section \ref{sec:GIFSdege}. We consider the following abstract setting for perturbed graph-directed systems:
\ite
{
\item[(D.1)] The unperturbed system $(G,(J_{v}),(O_{v}),(T_{e}))$ satisfies conditions $(P.1)$, $(G.5)_{C}$, $(G.7)_{B}$, $(G.8)_{SV}$, $(G.9)_{J}$, and either $(G.6)_{S}$ or $(G.6)_{C}$. 
\item[(D.2)] The perturbed system $(\hat{G},(J_{v}(\e))_{v\in \hat{V}},(O_{v}(\e))_{v\in \hat{V}},(T_{e}(\e,\cd))_{e\in \hat{E}})$, depending on a small parameter $\e\in (0,1)$, satisfies conditions $(P.2)$ and $(P.3)$. Moreover, for each $\e\in (0,1)$, the system satisfies $(G.8)_{SV}$, $(G.9)_{J}$, and either $(G.6)_{S}$ or both $(G.5)_{C}$ and $(G.6)_{O}$.
\item[(D.3)] The potential $\ph:=\ph^{V}$ is regular, that is, there exists $s_{0}\in \R$ such that $P(s_{0}\ph)=0$. Moreover, $s_{0}\in I_{P}^{\hat{V}}$, where  $I_{P}^{\hat{V}}$ is defined in (\ref{eq:I**=}).
\item[(D.4)] There exists $s\in \R$ such that $0<P(s\ph)<+\infty$ and $s\in \underline{I}_{P}^{\hat{V}}$, where $\underline{I}_{P}^{\hat{V}}$ is defined by (\ref{eq:uI**=}). In particular, $\ph$ is strongly regular.
}
Let $K(\e)$ denote the limit set of the perturbed GIFS. Then we obtain the following:
\thm
{\label{th:conv_dim}
Assume that conditions (D.1) and (D.2) are satisfied.
\ite
{
\item If condition (D.3) holds, then $\limsup_{\e\to 0}\dim_{H}K(\e)\leq \dim_{H}K$.
\item If condition (D.4) holds, then $\liminf_{\e\to 0}\dim_{H}K(\e)\geq \dim_{H}K$.
}
Consequently, if conditions (D.1)-(D.4) are satisfied, then $\lim_{\e\to 0}\dim_{H}K(\e)=\dim_{H}K$.
}
\pros
(1) Recall the potential $\phe$ defined in (\ref{eq:phe=...}). By condition (D.3), we have $s\in I_{P}^{\hat{V}}$ for any $s\geq s_{0}$. Fix any $\delta>0$. Then by Theorem \ref{th:conv_phe}, we have
$P((s_{0}+\delta)\phe)\to P((s_{0}+\delta)\ph)<0$.
Therefore, there exists $\e_{0}>0$ such that for any $0<\e<\e_{0}$, we have $P((s_{0}+\delta)\phe)<0$. Applying Theorem \ref{th:GIFS_upper_anygraph}, we obtain $\dim_{H}K(\e)\leq s(\e)<s_{0}+\delta$, proving the assertion.
\smallskip
\\
(2) Recall $\underline{\ph}(\e,\cd)$ defined in (\ref{eq:uphe=}). Let $s_{0}\in \R$ be such that $P(s_{0}\ph)=0$. By (D.4), there exists $\delta_{0}>0$ such that $s_{0}-\delta_{0}\in \underline{I}_{P}^{\hat{V}}$, and thus $s\in \underline{I}_{P}^{\hat{V}}$ for any $s\geq s_{0}-\delta_{0}$. Fix any $0<\delta<\delta_{0}$. By Theorem \ref{th:conv_phei}, we obtain
$P((s_{0}-\delta)\underline{\ph}(\e,\cd))\to P((s_{0}-\delta)\underline{\ph}(\e,\cd))<0$.
Therefore, there exists $\e_{0}>0$ such that for any $0<\e<\e_{0}$, we have $P((s_{0}-\delta)\underline{\ph}(\e,\cd))>0$. Then, Theorem \ref{th:GIFS_lower} gives $\dim_{H}K(\e)\geq s(\e)>s_{0}-\delta$. Hence the proof is complete.
\proe
\section{Applications}\label{sec:app}
\subsection{An example with a non-finitely irreducible transition matrix}\label{sec:ex_non-finitely}
We consider the following graph-directed system:
\ite
{
\item The graph $G=(V,E)$ is defined by $V=\{1,2,\dots\}$ and $E=\{1u\in V^{2}\,:\,u\geq 1\}\cup \{vu\in V^{2}\,:\,u=v-1,\ v\geq 2\}$.
\item For each $v\in V$, let $J_{v}=B(c_{v},1/2^{v})\subset \R^{D}$ and $O_{v}=B(c_{v},1)\subset \R^{D}$, where the centers $c_{v}$ are chosen so that $\mathrm{int}J_{v}\cap \mathrm{int}J_{v^\p}=\emptyset$ whenever $v\neq v^\p$.
\item For $vu\in E$, the map $T_{vu}\,:\,\R^{D}\to \R^{D}$ is a similarity with contraction ratio $(1/2)^{v}$, satisfying $T_{u}J_{u}\subset J_{v}$, and the OSC holds.
}
Then we obtain the following:
\thm
{\label{th:concex_non-finitely}
For the above GIFS $(G,(J_{v}),(O_{v}),(T_{e}))$, the Hausdorff dimension of its limit set $K$ is the unique solution $s$ to the equation $P(s\ph^{V})=0$. In particular, we have $\underline{I}^{V}=(0,\infty)$ and $\underline{I}^{E}=\emptyset$.
}
\pros
To apply Theorem \ref{th:GIFS_Bowen}, it suffices to verify condition $(G.8)_{SV}$. Note that $\underline{I}^{V}=\{s\geq 0\,:\,\sum_{k=1}^{\infty}(1/2^s)^{k}<+\infty\}=(0,\infty)$. To check the condition, consider the subgraph $G_{0}=(V_{0}=\{1,2\},E_{0}=\{11,12,21\})$ and the corresponding subsystem $(G_{0},(J_{v}),(O_{v}),(T_{e}))$. The limit set $K_{0}$ of this subsystem has the Hausdorff dimension $0.5514\cdots$, which is the unique solution $s$ of the equation $(1/2)^{s}+(1/8)^s=1$. Thus, $\dim_{H}K\geq \dim_{H}K_{0}>0$, and $(G.8)_{SV}$ is satisfied. Hence, the conclusion follows from Theorem \ref{th:GIFS_Bowen}.
\proe
\subsection{Perturbed affine maps with degeneration}\label{sec:PAMD}
We present an example of perturbed affine maps with degeneration:
\ite
{
\item[(I)] An unperturbed system $(G,(J_{v}),(O_{v}),(T_{e}))$ satisfies the following:
\ite
{
\item[(i)] A graph $G=(V,E,i,t)$ is simple and contain at least one strongly connected component.
\item[(ii)] For each $e\in E$, the map $T_{e}\,:\,\R^{D}\to \R^{D}$ is conformal and affine of the form $T_{e}(x)=M_{e}x+b_{e}$, where $M_{e}$ is a $D\times D$ matrix and $b_{e}\in \R^{D}$. Assume the following: $T_{e}J_{t(e)}\subset J_{i(e)}$, $\sup_{e\in E}\|M_{e}\|<1$, $\sum_{v\in V}\sup_{e\in E\,:\,i(e)=v}\|M_{e}\|^{s}<+\infty$ for any $s>0$, and for any $v\in V$, $\diam(J_{v})\leq c_{\adr{cjT}}\sup_{e\,:\,i(e)=v}\|M_{e}\|$ for some constant $c_{\adl{cjT}}>0$. Moreover, the system is strongly regular and satisfies the OSC.
}
\item[(II)] A perturbed system $(\hat{G},(J_{v}),(O_{v}),(T_{e}(\e,\cd)))$ satisfies the following:
\ite
{
\item[(i)] The graph $\hat{G}=(\hat{V},\hat{E},i,t)$ is simple and contains $G$ as a subgraph.
\item[(ii)] Each map $T_{e}(\e,\cd)\,:\,J_{t(e)}\to J_{i(e)}$ $(e\in \hat{E})$ have the form
\ali
{
T_{e}(\e,x)=
\case
{
(M_{e}+\ti{M}_{e}(\e))x+b_{e}(\e),&\text{if }e\in E\\
\ti{M}_{e}(\e)x+b_{e}(\e),&\text{if }e\in \hat{E}\setminus E,\\
}
}
where $\sup_{e\in \hat{E}}\|\ti{M}_{e}(\e)\|\to 0$, $\sum_{v\in \hat{V}}\sup_{e\,:\,i(e)=v}\sup_{\e>0}\|\ti{M}_{e}(\e)\|^{s}<+\infty$ for any $s>0$, there exist vectors $b_{e}\in \R^{D} (e\in \hat{E}\setminus E)$ such that $\sup_{e\in \hat{E}}|b_{e}(\e)-b_{e}|\to 0$. Moreover the SSC is satisfied for each $\e>0$.
}
}
Let $O_{v}=B(J_{v},\eta)$ for some fixed $\eta>0$. Since $T_{e}(\e,\cd)$ is a contraction uniformly in $e\in \hat{E}$ for sufficiently small $\e>0$, the inclusion $T_{e}(\e,O_{t(e)})\subset O_{i(e)}$ is satisfied. We then obtain the following:
\thm
{\label{th:conv_affine}
Assume that conditions (I) and (II) above are satisfied. Then we have $\dim_{H}K(\e)\to \dim_{H}K$ as $\e\to 0$.
}
\pros
It is straightforward to verify that conditions (D.1) and (D.2) of Theorem \ref{th:conv_dim} are valid. To check (D.3) and (D.4), we note that both $I_{P}^{\hat{V}}$ and $\underline{I}_{P}^{\hat{V}}$ are equal to $(0,+\infty)$. Indeed, for $0<s\leq 1$, we have $\sum_{v\in \hat{V}}\sup_{e\,:\,i(e)=v}\|M_{e}+\ti{M}_{e}(\e)\|^{s}\leq \sum_{v\in \hat{V}}\sup_{e\,:\,i(e)=v}\|M_{e}\|^{s}+\sum_{v\in \hat{V}}\sup_{e\,:\,i(e)=v}\|\ti{M}_{e}(\e)\|^{s}<+\infty$. Moreover, since $\ph$ is strongly regular, there exist $s_{0},\delta>0$ with $s_{0}-\delta>0$ such that $P(s_{0}\ph)=0$ and $0<P(s_{0}-\delta)<+\infty$. In addition to the fact $I_{P}^{\hat{V}}=\underline{I}_{P}^{\hat{V}}=(0,+\infty)$, conditions (D.3) and (D.4) are fulfilled. Hence the assertion follows from Theorem \ref{th:conv_dim}.
\proe
\subsection{Perturbed complex continued fractions with degeneration}\label{sec:PCCF}
In this section, we consider a GIFS driven by complex continued fractions, and formulate a perturbation of this system involving degeneration. Let $\hat{V}=\{v\}$ be a singleton vertex set, $\hat{E}\subset E_{*}:=\{m+n\sqrt{-1}\,:\,(m,n)\in \Z\times \Z,\ m\geq 1\}$ a nonempty countable edge set, $X_{v}=\overline{B(1/2,1/2)}$, and $O_{v}=B(1/2,3/4)$. For each $e\in \hat{E}$, we define a map $T_{e}\,:\,O_{v}\to O_{v}$ by $T_{e}(z)=1/(e+z)$.
Then the system $\hat{\mathcal{G}}=(\hat{G}=(\hat{V},\hat{E}),(J_{v}),(O_{v}),(T_{e}))$ is a conformal multi GIFS, with the exception that $T_{1}^\p(0)=1$ when $e=1\in \hat{E}$. However, this not pose a problem because the composition $T_{e}\circ T_{e^\p}$ is a contraction mapping uniformly in $ee^\p\in \hat{E}^{2}$, and therefore the condition in Remark \ref{eq:GIFS_rem} is satisfied. It is known in \cite[Section 7]{CLU2019} that when $\hat{E}=E_{*}$, we have $\inf \underline{I}_{P}^{\hat{V}}=1$, and the limit set of this system satisfies $\dim_{H}K\geq 1.825$ (see \cite{Priyadarshi2}). In particular, this system is strongly regular.

We now define perturbed maps $T_{e}(\e,z)$ under the following assumptions:
\ite
{
\item[(I)] A subgraph $G=(V,E)$ of the graph $\hat{G}=(\hat{V},\hat{E})$ satisfies that the subsystem $(G=(V,E),(J_{v})_{v\in V},(O_{v})_{v\in V},(T_{e})_{e\in E})$ associated with complex continued fractions has the limit set $K$ with positive Hausdorff dimension.
\item[(II)] For each $e\in \hat{E}$ and $\e\in (0,1)$, a function $T_{e}(\e,\cd)\,:\,J_{t(e)}\to J_{i(e)}$ is defined by
\ali
{
T_{e}(\e,z)=
\case
{
T_{e}(z), &\text{if }e\in E\\
\di\frac{1}{e+\frac{1}{2}+\e(z-\frac{1}{2})}, &\text{if }e\in \hat{E}\setminus E.\\
}
}
\item[(III)] If $s_{0}$ satisfies $P(s_{0}\ph^{E})=0$, then $s_{0}\in \{s\geq 0\,:\,\sum_{e\in \hat{E}}\|T_{e}^\p\|^{s}<+\infty\}$.
}
\begin{figure}[!htb]
\centering
\includegraphics[scale=0.25,angle=0]{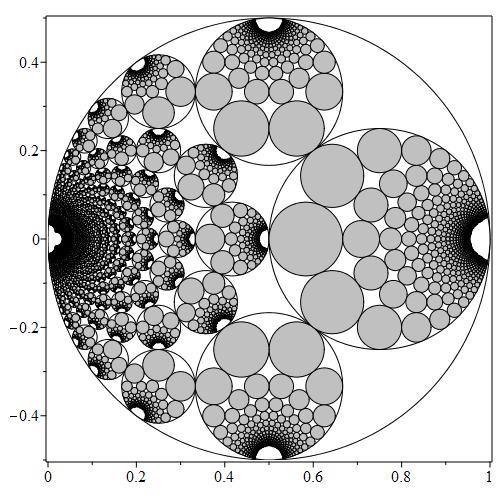}\qquad
\includegraphics[scale=0.25,angle=0]{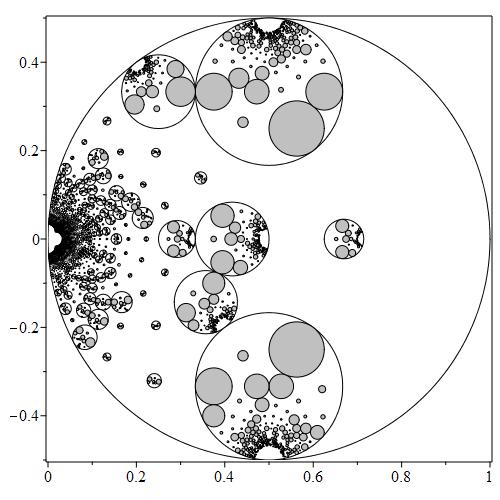}\qquad
\includegraphics[scale=0.25,angle=0]{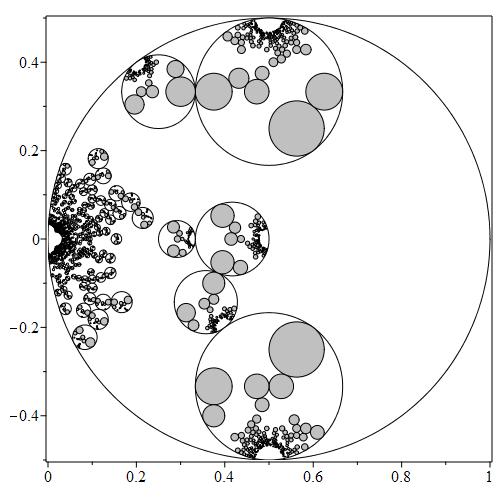}
\caption{Approximation of limit sets of the system of complex continued fractions (left), a perturbed subsystem (center), and the unperturbed subsystem (right).}\label{fig:2}
\end{figure}
Then we obtain the following:
\thms\label{th:per_CF_dege}
Assume that conditions (I)-(III) above are satisfied. Then the Hausdorff dimension of the limit set of the GIFS $(G,(J_{v}),(O_{v}),(T_{e}(\e,\cd)))$ converges to $\dim_{H}K$.
\thme
\pros
By Proposition \ref{prop:reduce_GIFS}, the unperturbed and perturbed systems can be reduced to simple GIFSs, denoted by $\tilde{\mathcal{G}}$ and $\tilde{\mathcal{G}}(\e)$, respectively. In particular, the new compact set $\tilde{J}_{v}(\e)=T_{e}(\e,J_{t(e)})$ satisfies condition (P.2)(ii), due to the uniform convergence of $T_{e}(\e,\cd)$. 
It is straightforward to verify that condition (D.1) in Theorem \ref{th:conv_dim} holds. To verify condition (D.2), observe that the new system $\tilde{\mathcal{G}}(\e)$ satisfies conditions $(G.5)_{C}$, $(G.6)_{O}$, $(G.7)_{B}$, and $(G.9)_{J}$ for each $\e\in (0,1)$. Moreover, since $0=P(s_{0}\ph^{V})\leq P(s_{0}\ph(\e,\cd))<+\infty$, and condition (III) ensures that $s_{0}$ belongs to the summability domain, condition $(G.8)_{SV}$ is also satisfied.
Note that $T_{e}(\e,z)$ converges to the point $1/(e+1/2)$ uniformly in $z\in J_{t(e)}$ and $e\in \hat{E}\setminus E$. From the identity $T_{e}(\e,z)=T_{e}((1/2)+\e(z-1/2))$, we have 
\alil
{
|T_{e}^\p(\e,z)|=\e|T_{e}^\p((1/2)+\e(z-1/2))|=\frac{\e}{|e+(1/2)+\e(z-1/2)|^{2}}.\label{eq:diff_Tepez}
}
Then, it follows that $|T_{e}^\p(\e,z)|$ vanishes uniformly in $z$ and $e\in \hat{E}\setminus E$. The bounded distortion $\||T_{e}^\p(\e,x)|-|T_{e}^\p(\e,x)|\leq c_{\adr{bd3}}|T_{e}^\p(\e,x)||x-y|^{\beta}$ uniformly in $e\in \hat{E}$ are also satisfied. Therefore conditions (P.2) and (P.3) are satisfied, and thus (D.2) is verified. To verify condition $(D.3)$, we note that (\ref{eq:diff_Tepez}) yields
\alil
{
I:=I_{P}^{\hat{V}}=\underline{I}_{P}^{\hat{V}}=\{s\geq 0\,:\,\sum_{e\in \hat{E}}\|T_{e}(\e,\cd)^\p\|^{s}<+\infty\}=&\{s\geq 0\,:\,\sum_{e\in \hat{E}}\|T_{e}^\p\|^{s}<+\infty\}\label{eq:IVIV=}\\
=&\{s\geq 0\,:\,P(s\phe)<+\infty\}\nonumber
}
for any $\e\in (0,1)$, where the fourth equation is implied by (\ref{eq:diff_Tepez}), and the last equation holds from the transition matrix of $\hat{G}$ is finitely irreducible. Hence (D.3) is satisfied. From Theorem \ref{th:conv_dim}(1), it follows that $\limsup_{\e\to 0}\dim_{H}K(\e)\leq \dim_{H}K$. On the other hand, if $\inf I<\dim_{H}K$, then condition (D.4) holds, and therefore $\lim_{\e\to 0}\dim_{H}K(\e)=\dim_{H}K$. Finally, in the case $\inf I=\dim_{H}K$, for any $\delta>0$, we have $P((\dim_{H}K-\delta)\phe)=+\infty$, implying $\dim_{H}K-\delta\leq \dim_{H}K(\e)$. Letting $\delta\to 0$, we conclude $\dim_{H}K\leq \dim_{H}K(\e)$, and hence the convergence $\dim_{H}K(\e)\to \dim_{H}K$ holds.
\proe
If condition (III) does not hold, then there exists an example in which the Hausdorff dimension of $\dim_{H}K(\e)$ is discontinuous at $\e=0$. Note that condition (III) is always satisfied in the case of finite graphs, and hence discontinuity of the Hausdorff dimension does not occur in that setting.
\prop
{\label{prop:per_CF_dege_discon}
Assume that (I) and (II) are satisfied. Then there exists an example for which $\lim_{\e\to 0}\dim_{H}K(\e)\neq \dim_{H}K$.
}
\pros
For instance, consider the case where $\hat{E}=E_{*}$ and $0<\dim_{H}K<1$ (this case exists by \cite[Theorem 1.4]{CLU2019}). Then, by $\hat{E}=E_{*}$ and (\ref{eq:IVIV=}), we obtain $\inf I=1$. Hence $\dim_{H}K(\e)\geq 1$ for any $\e>0$, and consequently $\dim_{H}K(\e)$ does not converge to $\dim_{H}K<1$.
\proe

\end{document}